\documentclass[a4paper,12pt]{amsart}
\usepackage[utf8]{inputenc}
\usepackage[T1]{fontenc}
\usepackage[UKenglish]{babel}
\usepackage[a4paper,margin=28mm]{geometry}
\usepackage{verbatim}

\allowdisplaybreaks[4]
\usepackage{times}
\usepackage{dsfont,mathrsfs}
\usepackage{amsmath}
\usepackage{amsthm}
\usepackage{amssymb}
\usepackage{amsfonts}
\usepackage{latexsym}
\usepackage{booktabs}

\usepackage[colorlinks,linkcolor=red,pagebackref]{hyperref}
\usepackage{xcolor}

\usepackage{ulem}

\newtheorem{theorem}{Theorem}[section]
\newtheorem{lemma}[theorem]{Lemma}
\newtheorem{proposition}[theorem]{Proposition}
\newtheorem{corollary}[theorem]{Corollary}
\newtheorem{assumption}[theorem]{Assumption}

\theoremstyle{definition}
\newtheorem{definition}[theorem]{Definition}

\newtheorem{remark}[theorem]{Remark}
\numberwithin{equation}{section}
\renewcommand{\labelenumi}{\roman{enumi})}
\renewcommand\theenumi\labelenumi
\renewcommand{\leq}{\leqslant}

\renewcommand{\geq}{\geqslant}

\newcommand{\tl}{\tilde}
\newcommand{\Be}{\begin{equation}}
\newcommand{\Ees}{\end{equation*}}
\newcommand{\Bes}{\begin{equation*}}
\newcommand{\Ee}{\end{equation}}

\newcommand{\R}{\mathbb{R}}
\newcommand{\Rd}{\mathbb{R}^d}
\newcommand{\E}{\mathbb{E}}

\newcommand{\e}{\varepsilon}
\newcommand{\ep}{\epsilon}
\newcommand{\PP}{\mathbb{P}}

\newcommand{\Gcal}{\mathcal{G}}
\newcommand{\Lcal}{\mathcal{L}}

\newcommand{\dif}{\mathrm{d}}

\usepackage{marginnote}
\begin{document}
\title[Non-asymptotic Analysis of Poisson randomized midpoint Langevin Monte Carlo]
{Non-asymptotic Analysis of Poisson randomized midpoint Langevin Monte Carlo}

\author[T. Shen]{Tian Shen}
\address{Tian Shen: School of Statistics and Data Science, Shanghai University of International Business and Economics, Shanghai, China;}
\email{shentian@zju.edu.cn}

\author[Z. Su]{Zhonggen Su*}
\address{Zhonggen Su: School of Mathematical Sciences, Zhejiang University, Hangzhou, Zhejiang, China;}
\email{suzhonggen@zju.edu.cn}


\begin{abstract}

The task of sampling from a high-dimensional distribution $\pi$ on $\R^d$ is  a  fundamental algorithmic problem with applications throughout statistics, engineering, and the sciences. 
Consider the Langevin diffusion on $\R^d$
\begin{align*}
\dif X_t=-\nabla U(X_t)dt+\sqrt{2}dB_t,
\end{align*}
under mild conditions, it admits $\pi(\dif x)\propto \exp(-U(x))\dif x$ as its unique stationary distribution. Recently, Kandasamy and   Nagaraj (2024) introduced  a   stochastic algorithm  called Poisson Randomized Midpoint  Langevin Monte Carlo (PRLMC)  to enhance the rate of convergence towards the  target distribution $\pi$. In this paper, we first show that under mild conditions, the PRLMC, as a Markov chain, admits a unique stationary distribution $\pi_\eta$ ($\eta$ is the step size) and   obtain the convergence rate of PRLMC to $\pi_\eta$ in total variation distance. Then we establish an sharp error bound between $\pi_\eta$ and $\pi$ under the 2-Wasserstein distance. Finally, we propose a  decreasing-step size version of PRLMC  and provide its convergence rate to $\pi$ which is nearly optimal.
\end{abstract}

\maketitle


\section{Introduction} \label{Model}
Our primary aim is to study convergence properties of Markov chains that arises as approximations of ergodic solutions of the Langevin diffusion
\begin{align}\label{ld}
\dif X_t=-\nabla U(X_t)dt+\sqrt{2}dB_t, \quad X_0\in\R^d,
\end{align}
where $U:\R^d\to\R$ is the potential function and $(B_t)_{t\geq0}$ is a $d$-dimensional Brownian motion. It is well-known that, under some mild conditions, \eqref{ld} admits a strong solution $(X_t^x)_{t\geq0}$ for any initial state $x\in\R^d$, see Roberts and Tweedie \cite{ro}. Moreover, $(X_t^x)_{t\geq0}$ admits $\pi(\dif x)=\frac{1}{Z} \exp(-U(x))\dif x$ with normalized constant $Z=\int_{\R^d}\exp(-U(x))\dif x$ as its unique stationary distribution. Sampling from such $\pi$ is a fundamental problem which occurs for example in Bayesian inference and machine learning, see Cotter et al. \cite{an} and Cotter et al. \cite{co} for details. The most commonly used method in literature is the Euler-Maruyama discretization:
\begin{align}\label{em}
X_{(k+1)\eta}=X_{k\eta}-\eta\nabla U(X_{k\eta})+\sqrt{2\eta}\xi_{k+1},
\end{align}
where $\eta$ is the step size and $\{\xi_k\}_{k\in\mathbb N}$ is a sequence of i.i.d. $d$-dimensional standard Gaussian random vectors. Following Roberts and Tweedie \cite{ro}, this recursion is called the Unadjusted Langevin Algorithm (ULA). Under mild conditions, the Markov chain $(X_{k\eta})_{k\in\mathbb N_0}$ defined in \eqref{em} admits a unique stationary distribution $\pi_\eta$. 
Although $\pi_\eta$ generally differs from $\pi$, they are close   when the step size $\eta$ is sufficiently small. Numerous articles have studied the  nonasymptotic bounds between the two stationary distributions. For instance,  Dalalyan \cite{da}, Durmus and Moulines \cite{du, du2},  and Fang et al.\cite{fa} established such bounds under total variation distance, weighted total variation distance, 2-Wasserstein distance and 1-Wasserstein distance, respectively.

In practical applications, people often consider the issue of algorithmic complexity, which involves determining how many iterations the Markov chain \eqref{em} requires so that its distribution is within a given error tolerance $\epsilon$ of the target $\pi$. This problem has been extensively studied in literatures (such as \cite{da2,da3,du00,ve}) under various assumptions on the potential function $U$ such as log-concavity, Log-Sobolev inequality and isoperimetry. To enhance the convergence rate, Shen and Lee \cite{sh} proposed the randomized Langevin monte Carlo (RLMC):
\begin{align*}
X_{(k+1)\eta}&=X_{k\eta}-\eta\nabla U(X_{(k+u_k)\eta})+\sqrt{2\eta}\xi_{k+1},
\end{align*}
where
\begin{align*}
X_{(k+u_k)\eta}&=X_{k\eta}-u_k\eta\nabla U(X_{k\eta})+\sqrt{2u_k\eta}\xi_{k+1}^\prime,
\end{align*}
and $(u_k)_{k\in\mathbb N_0}$ is a sequence of i.i.d. uniform random variables on the interval $[0,1]$. Consider  the Langevin diffusion from time $k\eta$ to $(k+1)\eta$,
\begin{align}\label{1334}
X_{(k+1)\eta}=X_{k\eta}-\int_{k\eta}^{(k+1)\eta}\nabla U(X_t)\dif t+\sqrt2(B_{(k+1)\eta}-B_{k\eta}).
\end{align}
The difference between ULA and RLMC lies in the fact that ULA uses $\eta\nabla U(X_{k\eta})$ to approximate the integral term in RHS of \eqref{1334}, whereas RLMC uses $\eta\nabla U(X_{(k+u_k)\eta})$, which is a unbiased estimator of the integral.
The  randomized midpoint  discretization  is  developed further in \cite{al,al2,hy,yl}, and it is showed that in some sense the randomized midpoint discretization  converges to target $\pi$ faster than the general ULA.

Recently, based on the randomized midpoint discretization,  Kandasamy and   Nagaraj \cite{ka} proposed a  stochastic algorithm called Poisson Randomized Midpoint  Langevin Monte Carlo (PRLMC): given step size $\eta$ and a positive integer $K$,  
\begin{align}\label{prlmc}
\tl{X}_{k+1}=\tl{X}_{k}-\eta\nabla U(\tl{X}_{k})+\eta\sum_{i=0}^{K-1}H_{k,i}(\nabla U(\tl{X}_{k})-\nabla U(\hat{X}_{k,i}))+\sqrt{2\eta}\xi_{k+1}
\end{align}
where $\{H_{k,i}\}_{k\in\mathbb N_0, 0\leq i\leq K-1}$ is a sequence of i.i.d. Bernoulli random variables with parameter $1/K$ and  for $i=0,1,\ldots,K-1$, 
\begin{align*}
\hat{X}_{k,i}=\tl{X}_{k}-\frac{i\eta}{K}\nabla U(\tl{X}_{k})+\sqrt{\frac{2i\eta}{K}}\xi_i.
\end{align*}
Kandasamy and   Nagaraj \cite{ka} established an error bound in KL divergence  between the law of PRLMC and that of Langevin diffusion. Srinivasan and Nagaraj \cite{sr} gave an  error bound between ${\rm law}(\tl X_k)$ and ${\rm law}(X_k)$ in terms of 2-Wasserstein distance, where $\tl X_k$ and $X_k$ represents the $k$-th iteration in PRLMC and ULA respectively. And the authors also derived an error bound between ${\rm law}(\tl X_k)$ and the target $\pi$ based on this result. 

Our focus in this paper is on the PRLMC $(\tl X_k)_{k\in\mathbb N_0}$ defined in \eqref{prlmc}. It is called Poisson midpoint method since the number of midpoints in any interval $[k\eta, (k+1)\eta]$ converges to a Poisson distribution with parameter 1 as $K\to\infty$, i.e. for  nonnegative integer $n$
\begin{align*}
&\lim_{K\to\infty}\PP\Big({\rm number \ of \  midpoints \ in \ interval \ } [k\eta,(k+1)\eta]=n\Big)\\
=&\lim_{K\to\infty}\binom{K}{n}(\frac{1}{K})^n(1-\frac{1}{K})^{K-n}\\
=&\frac{1}{n!e}.
\end{align*}

   Given $\tl{X}_{0}=x$, then $\hat{X}_0=x$ and  for $i=1,2,\ldots,K-1$, $\hat{X}_i$ (we omit the subscript 0 for convenience) is Gaussian with density
\begin{align*}
\rho_i(y_i)=&\PP(\hat{X}_i=y_i|\tl{X}_{0}=x)\\
=&\frac{1}{(\sqrt{2\pi})^d(\frac{2i\eta}{K})^{d/2}}\exp\left(-\frac{[y_i-(x-\frac{i\eta}{K}\nabla U(x))]^2}{\frac{4i\eta}{K}}\right)
\end{align*}

According to \eqref{prlmc}, $(\tl{X}_k)_{k\geq0}$ is a Markov chain on  $\R^d$. Given the values of $\{H_i\}_{i=0}^{K-1}$,
\begin{align*}
&q_\eta(x,\tl x)\big |_{H_0,H_1,\cdots,H_{K-1}}\\
=&\int_{(\R^{d})^{K-1}}
\frac{1}{(4\pi\eta)^{d/2}}\exp\left(-\frac{\left[\tl x-\left(x-\eta\nabla U(x)+\eta\sum\limits_{i=1}^{K-1}H_i(\nabla U(x)-\nabla U(y_i))\right)\right]^2}{4\eta}\right)\\
&\qquad \cdot\Pi_{i=1}^{K-1}\rho_i(y_i)\dif y_1\cdots\dif y_{K-1}.
\end{align*}
 
So the one-step transition probability density of $(\tl{X}_k)_{k\in\mathbb N_0}$ is given by
\begin{align}\label{zymd}
q_\eta(x,\tl x)=\E\left[q(x,\tl x) |_{H_0,H_1,\cdots,H_{K-1}}\big |H_0,H_1,\cdots,H_{K-1}\right].
\end{align}
and transition kernel 
\begin{align}\label{zyh}
Q_\eta(x,A)=\int_{A}q_{\eta}(x,\tl x)\dif\tl x,
\end{align}
for any $x\in\R^d$ and $A\in\mathcal B(\R^d)$.

The contributions of this work primarily include:

(1) We prove that, under mild conditions, $(\tl X_n)_{n\in\mathbb N_0}$ admits a unique stationary distribution $\pi_\eta$. Furthermore, we  show the convergence rate of the marginal distribution of this Markov chain to $\pi_\eta$ in total variation distance is $d_{\rm TV}(\mathcal L(\tl X_n),\pi_\eta)=o(\delta^{-n}) , \quad n\to\infty$;

(2) We derive an upper bound  $ O(\sqrt\eta)$ of the  2-Wasserstein distance between $\pi_\eta$ and  $\pi$ under Assumption \ref{ass};  moreover, we obtain a sharp bound $O(\eta)$ if additionally $\|\nabla^3U\|<\infty$ ;

(3) We also design a decreasing-step size PRLMC to approximate $\pi$ and provide the convergence rates under 2-Wasserstein distance and a certain functional class distance.

\textbf{Notations} We end this section by introducing some notations, which will be frequently used in the sequel. The inner product of $x, y \in\R^d$ is denoted by $\langle x,y\rangle$ and the Euclidean metric is denoted by $|x|$.
For any $u\in\mathbb{R}^{d}$, matrix $A=(A_{ij})_{d\times d}$,  define $Au^{\otimes 2}=\langle Au, u\rangle$, and the operator norm of the matrix $A$ is denoted by
$\|A\|_{\rm op}=\sup_{|u|=1}|Au|$. Denote by $\vec{\Delta} F$ the vector Laplacian of $F$ defined by: for all $x\in\Rd$, $\vec{\Delta} F(x)$ is the vector of $\R^n$ such that for all $i\in\{1,\ldots,n\}$, the $i$-th component of $\vec{\Delta} F(x)$ equals to $\sum_{j=1}^d
(\partial^2F_i/\partial x_j^2)(x)$.

 Let $\mathcal{C}^k(\R^d,\R)$ with $k\geq 1$ denote the collection of all  $k$-th order continuously differentiable functions from $\R^d$ to $\R$.  For $f\in \mathcal{C}^k(\R^d,\R)$, the  supremum norm is defined as
\begin{align*}
\|\nabla^if\|=\sup_{x\in\R^d}\|\nabla^if(x)\|_{\rm op}, \quad  i=1,2.
\end{align*}

Let $(\Omega, \mathcal F, \PP)$ be the probability space, for any random vector $X$, the law of $X$ is denoted by $\mathcal L(X)$ and the $L^2$-norm is defined by $\|X\|_2=[\E(|X|^2)]^{\frac{1}{2}}$. Let $\mathcal P(\R^d)$ be the space of probability measures on $\R^d$. For $\mu,\nu\in\mathcal P(\R^d)$, their total variation distance is defined as
\begin{align*}
\dif_{\rm TV}(\mu,\nu)=\sup_{A\in\mathcal B(\R^d)}|\mu(A)-\nu(A)|,
\end{align*}
and the 2-Wasserstein distance is defined as
\begin{align*}
W_2(\mu,\nu)=\left(\inf_{\lambda\in\Pi(\mu,\nu)}\int_{\R^d\times\R^d}|x-y|^2\lambda(\dif x,\dif y)\right)^{\frac{1}{2}},
\end{align*}
where $\Pi(\mu,\nu)$ is the space of probability measures $\lambda\in\mathcal P(\R^d\times\R^d)$ whose marginals are $\mu$ and $\nu$ respectively.

Denote  a functional class
\begin{align*}
\mathcal G=\{h\in\mathcal{C}^2(\R^d,\R):\|\nabla^ih\|<+\infty \ {\rm with}\  i=1,2\}.
\end{align*}
Then for any $h\in\Gcal$, there exists a  positive constant $C_h$ such that $|h(x)|\leq C_h(1+|x|^2)$. Denote two probability measure spaces
\begin{align*}
\mathcal P_0=\left\{\mu\in\mathcal P(\R^d):\int_{\R^d}h(x)\mu(\dif x)<+\infty \  {\rm for\ any }\ h\in\Gcal\right\}
\end{align*}
and 
\begin{align*}
\mathcal P_2=\left\{\mu\in\mathcal P(\R^d):\int_{\R^d}|x|^2\mu(\dif x)<+\infty\right\},
\end{align*}
then it is easy to see that $\mathcal P_2\subset\mathcal P_0$. For $\mu_1,\mu_2\in\mathcal P_0$, define
\begin{align*}
\dif_{\Gcal}(\mu_1,\mu_2)=\sup_{h\in\Gcal}\left\{\int_{\R^d}h(x)\mu_1(\dif x)-\int_{\R^d}h(x)\mu_2(\dif x)\right\}.
\end{align*}
Obviously, $\dif_{\Gcal}(\cdot,\cdot)$ is a metric on $\mathcal P_0$.

The superscript $x$ in $X_t^x$ stands for the initial state.

The rest of paper is organized as follows. In Section \ref{mr}, we list some assumptions about potential function $U$ and step sequence $\gamma_n$,  and state our main results on the convergence rate of PRLMC. In Section \ref{1659},  we focus on the PRLMC with constant step size and prove Theorem \ref{ms} and Proposition \ref{ms2}. Then we turn to the decreasing step case and prove Theorems \ref{ms3} and \ref{ms4} in Sections \ref{2127} and \ref{sec5}, respectively.
 Appendix is  devoted to proofs of  auxiliary results stated in the previous Sections.

\section{Assumptions and Main Results}\label{mr}
 Consider the following assumption on the potential $U$.
\begin{assumption}\label{ass}
(i)  $U$ is strongly convex and gradient Lipschitz:
\begin{align}\label{m}
\frac{m}{2}|x-y|^2\leq U(y)-U(x)-\langle\nabla U(x),y-x\rangle\leq\frac{L}{2}|x-y|^2,\quad x,y\in\R^d
\end{align}

(ii) 
$\nabla U(0)=0$.

\end{assumption}
\begin{remark}
It follows from \cite[Lemma 4, Lemma 5]{dw} that \eqref{m} is equivalent to
\begin{align*}
|\nabla U(x)-\nabla U(y)|\leq L|x-y|
\end{align*}
and 
\begin{align*}
\langle -\nabla U(x)+\nabla U(y),x-y\rangle \leq -m|x-y|^2.
\end{align*}
Taking $y=0$ implies that for any $x\in\R^d$,
\begin{align}\label{u1}
|\nabla U(x)|^2\leq L^2|x|^2
\end{align}
and
\begin{align}\label{u2}
\langle -\nabla U(x),x\rangle \leq -m|x|^2.
\end{align}
 Assume $U$ is strongly convex, \cite[Theorem 2.1.8]{ne} shows that $U$ has a unique minimizer $x^*\in\R^d$ satisfying $\nabla U(x^*)=0$. In  Assumption \ref{ass} (ii),   we assume $x^*=0$, which is merely for computational convenience; $\nabla U(0)\neq 0$   does not affect the main results of this paper.
\end{remark}
Throughout this paper, we set 
\begin{align}\label{kappa}
    \kappa=\frac{2mL}{m+L}.
\end{align}

Our first result is the following theorem on the existence and uniqueness of the invariant distribution  of Markov chain $\tl{X}_{n}$ \eqref{prlmc}, and we also obtain the convergence  rate of $\mathcal L(\tl X_n)$ to its stationary distribution in total variation distance.

\begin{theorem}\label{ms}
Under Assumption \ref{ass}, there exists a constant $\eta_0\in(0,1)$ depending only on $m$ and $L$, such that for any $\eta\in(0,\eta_0)$, the Markov chain  $(\tl X_n)_{n\geq0}$ defined in \eqref{prlmc} has a unique stationary distribution $\pi_{\eta}$  which satisfies 
\begin{align}\label{1149}
\pi_\eta(|x|^2)=:\int_{\R^d}|x|^2\pi_\eta(\dif x)\leq\frac{d}{m}\{2+4\eta^2L^2(2-1/K)+L^2(1-1/K)\}.
\end{align}

 Furthermore, for any initial state $\tl X_0=x\in\R^d$, we have
\begin{align*}
\sum_{n=1}^\infty\delta^nd_{\rm TV}(\mathcal L(\tl X_n),\pi_\eta)
<\infty,
\end{align*}
where $\delta>1$ is a constant depending on step size $\eta$, and then the convergence rate of $(\tl X_n)_{n\geq0}$ to its stationary distribution is 
\begin{align*}
d_{\rm TV}(\mathcal L(\tl X_n),\pi_\eta)=o(\delta^{-n}), \quad n\to\infty.
\end{align*}
\end{theorem}

A natural  question arises: how far is $\pi_\eta$ away from $\pi$? Towards that, we derive an upper bound of 2-Wasserstein distance between $\pi_\eta$ and $\pi$.

\begin{proposition}\label{ms2}
Under the above notations and Assumption \ref{ass} ,we have that for any $\eta\in(0,\eta_0\wedge \frac{2}{m+L})$,
\begin{align*}
W_2(\pi_\eta,\pi)\leq\frac{L}{m}\eta^{\frac{1}{2}}\left[(2d+4\eta^2Ld+8L^2\eta^2d)^{\frac{1}{2}}+\sqrt{4-2/K} \left(\eta L^2\pi_\eta(|x|^2)/3+ d\right)^{\frac{1}{2}}\right].
\end{align*}
\end{proposition}

\begin{remark}
Theorem \ref{ms} and Proposition \ref{ms2} show that the law of  PRLMC \eqref{prlmc}  converges to $\pi_\eta$, not the real target distribution $\pi$. 
Plugging the bound of $\pi_{\eta}(|x|^2)$  from  \eqref{1149} into the above inequality shows that  the bias between the stationary distribution of the Langevin diffusion and that of the PRLMC is of the order $O(\sqrt{\eta})$.
\end{remark}

Next, we introduce the Poisson randomized midpoint  Euler-Maruyama scheme  with decreasing step-size sequence $\Gamma=(\gamma_n)_{n\in \mathbb{N}}$ below.
Let  $\gamma_k$ be the $k$-th step size, $t_0=0$, $t_k=\sum_{i=1}^k \gamma_i$, $k \geq 1$. Denote
\begin{align}\label{bbc}
Y_{t_{k+1}}=Y_{t_{k}}-\gamma_{k+1}\nabla U(Y_{t_{k}})+\gamma_{k+1}\sum_{i=0}^{K-1}H_{k,i}(\nabla U(Y_{t_{k}})-\nabla U(\hat Y_{t_{k},i}))+\sqrt{2\gamma_{k+1}}\xi_{k+1}.
\end{align}
Here $\{H_{k,i}\}_{k\in\mathbb N, 0\leq i\leq K-1}$  is a sequence of independent Bernoulli random variables with parameter $1/K$, $(\xi_k)_{k\geq1}$ is a sequence of i.i.d. $d$-dimensional standard normal vectors and  for $i=0,1,\ldots,K-1$, 
\begin{align*}
\hat Y_{t_{k},i}=Y_{t_{k}}-\frac{i\gamma_{k+1}}{K}\nabla U(Y_{t_{k}})+\sqrt{\frac{2i\gamma_{k+1}}{K}}\xi_i.
\end{align*}

\begin{assumption}\label{assump2}
(i) Let $(\gamma_n)_{n\in \mathbb{N}}$ be non-increasing and satisfy  $\gamma_n>0, \forall n\in \mathbb{N}$,  $\lim\limits_{n\to\infty}\gamma_n=0$ and  $\sum\limits_{n\geq1}\gamma_n=+\infty;$

(ii)
\begin{align*}
\omega:=\limsup_{n\rightarrow\infty}\frac{\gamma_{n}-\gamma_{n+1}}{\gamma_{n+1}^2}<+\infty.
\end{align*}
\end{assumption}

\begin{assumption}\label{ass2}
The potential $U$ is three times continuously differentiable and there exists $\tl L$ such that for all  $\|\nabla^3U\|\leq \tl L$.
\end{assumption}

Note that under Assumption  \ref{ass2}, we have that for all $x,y\in\R^d$,
\begin{align}\label{1352sf}
 \|\nabla^2U(x)-\nabla^2U(y)\|_{\rm op}\leq \tl L |x-y|, \quad |\vec{\Delta}(\nabla U)(x)|^2\leq d^2{\tl L}^2.
\end{align}

\begin{theorem}\label{ms3}

Let the potential  $U$  satisfy Assumptions \ref{ass} and \ref{ass2}, and  $(X_t^x)_{t\geq0}$ be the Langevin dynamics with unique  stationary distribution $\pi$. Let the step sequence $\Gamma$ satisfy Assumption \ref{assump2} with $2\omega<m$  and $m$  in \eqref{m}, then we have
\begin{align}\label{convergence1}
d_{\mathcal{G}} \left(\mathcal L\big( X_{t_n}^{x}\big),\mathcal L\big( Y_{t_n}^{x} \big)\right) \leq C (1+|x|^2)\gamma_n.
\end{align}
Furthermore, there exists some positive constant $C$ independent of $\Gamma$ such that
\begin{align}\label{convergence2}
d_{\mathcal{G}} \left(\mathcal L\big(Y_{t_n}^{x} \big), \pi\right) \leq C (1+|x|^2) \gamma_n.
\end{align}
\end{theorem}

\begin{remark}
According to Lemma \ref{LYexp} and \eqref{1234} below, for any $n\geq1$, $\Lcal\big( X_{t_n}^{x}\big), \Lcal\big( Y_{t_n}^{x} \big)$ and $\pi$ belong to $\mathcal P_2$. So the left-hand side of  \eqref{convergence1} and \eqref{convergence2} are well-defined.
\end{remark}

Based  on the idea of  Durmus and Moulines \cite[Theorem 8]{du2},  we can obtain the nonasymptotic bound in 2-Wasserstein distance between $\mathcal L\big(Y_{t_n}^{x} \big)$ and $\pi $.

\begin{theorem}\label{ms4}
Keep the same notations and assumptions  as in Theorem \ref{ms3}. Let $\Gamma=(\gamma_k)_{k\geq1}$ be a non-increasing step sequence with $\gamma_1\leq 1/(m+L)$.
 Then for any $n\geq1$,
\begin{align*}
W_2^2(\Lcal(Y_{t_n}^x),\pi)\leq u_n^{(1)}(\Gamma)\{|x|^2+d/m\}+u_n^{(2)}(\Gamma),
\end{align*}
 where $u_n^{(1)}(\Gamma)=2\prod\limits_{k=1}^{n}(1-\kappa\gamma_k/2)$ and
\begin{align}\label{10111201}
   u_n^{(2)}(\Gamma)=&\sum_{i=1}^{n}\Bigg[\gamma_i^3\Big\{L^2d[10-4/K+L^2\gamma_{i}^2/6+2\kappa^{-1}L^2\gamma_{i}]+6\kappa^{-1}d^2\tl L^2\nonumber\\
    &\qquad+ C(1+|x|^2)L^4[(8-4/K)\gamma_i+3\kappa^{-1}]\nonumber\\
&\qquad\left.+dL^4(\gamma_i+3\kappa^{-1})/m\Big\}\prod_{k=i+1}^n(1-\kappa\gamma_k/2)\right],
\end{align}
$\kappa$ is defined in \eqref{kappa} and the constant $C$ in $ u_n^{(2)}(\Gamma)$ is precisely the one appearing in \eqref{Yexp}.
\end{theorem}

\begin{corollary}\label{ms5}
Assume that  Assumptions \ref{ass} and \ref{ass2} hold. Let $\Gamma$ be a non-increasing step sequence with $\gamma_1\leq 1/(m+L)$, and suppose $\Gamma$ satisfies Assumption \ref{assump2} (i). Then $\lim_{n\to\infty}W_2(\Lcal(Y_{t_n}^x),\pi)=0.$ 
\end{corollary}

If $\gamma_k=\eta$ for all $k\geq1$, we can deduce from Theorem \ref{ms4}, a sharp bound between $\pi$ and the stationary distribution $\pi_\eta$ of PRLMC.

\begin{corollary}\label{ms6}
Under Assumptions \ref{ass} and \ref{ass2}, for any $\eta\in(0,\eta_0\wedge \frac{1}{m+L})$, we have
\begin{align*}
W_2^2(\pi,\pi_\eta)\leq&2\kappa^{-1}d\eta^2\Big\{L^2[10-4/K+L^2\eta^2/6+2\kappa^{-1}L^2\eta]+6\kappa^{-1}d\tl L^2\\
    &\qquad\quad \quad+ \{2+4\eta^2L^2(2-1/K)+L^2(1-1/K)\}/m\\
    &\qquad \quad \quad \ \ \cdot L^4[(8-4/K)\eta+3\kappa^{-1}]\\
    &\qquad\qquad+L^4(\eta+3\kappa^{-1})/m\Big\}.
\end{align*}
\end{corollary}

\section{ Proofs of Theorem \ref{ms} and Proposition \ref{ms2}}\label{1659}
\subsection{Preliminaries}
First we enumerate several definitions of Markon chain, which can be found in \cite{MC}. These definitions play a significant role in the proofs of our main results.  

\begin{definition}[Accessible set, small set]\label{small}
Let $P$ be a Markov transition kernel on state space $(X,\mathcal X)$.
\

(i) A set $A\in \mathcal X$ is said to be accessible if for all $x\in X$, there exists an integer $n\geq1$ such that $P^n(x,A)>0$.

(ii) A set $A\in \mathcal X$ is called a small set if there exist a positive integer $n$ and a nonzero measure $\mu$ on $(X,\mathcal X)$ such that for all $x\in A$ and $B\in\mathcal X$,
\begin{align*}
P^n(x,B)\geq\mu(B).
\end{align*}
The set $A$ is then said to be an $(n,\mu)$-small set.
\end{definition}

\begin{definition}[Irreducible]\label{irre}
A Markov kernel $P$ on $X\times\mathcal X$ is said to be irreducible if it admits an accessible small set.
\end{definition}

\begin{definition}[Period of an accessible small set]
The period of an accessible small set $C$ is the positive integer $d(C)$ defined by
\begin{align*}
d(C)=g.c.d.\{n\in\mathbb N^*:\inf_{x\in C}P^n(x,C)>0\},
\end{align*}
where $g.c.d.\{\cdot\}$ means the greatest common divisor of set $\{\cdot\}$.
\end{definition}

\begin{definition}[Period, aperiodicity, strong aperiodicity]\label{ape}
Let $P$ be an irreducible Markov kernel on $X\times\mathcal X$.

(i) The common period of all accessible small sets is called the period of $P$;

(ii) If the period is 1, $P$ is said to be aperiodic;

(iii) If there exists an accessible $(1,\mu)$-small set $C$ with $\mu(C)>0$, $P$ is said to be strongly aperiodic.
\end{definition}

\subsection{Auxiliary results}
For the Markov chain $(\tl X_k)_{k\geq0}$ defined in \eqref{prlmc}, we have the following properties which are similar to  Ye and Fan \cite[Proposition 3.1]{yi}.
\begin{proposition}\label{sms}
(i) If $E$ is a compact subset of $\R^d$ such that  $Leb(E)>0$, then $E$ is an accessible $(1,\epsilon\nu)$- small set, where
\begin{align*}
\epsilon=\inf_{x\in E,\tl x\in E}q_\eta(x,\tl x),
\end{align*}
$q_\eta$ is the transition density function defined in \eqref{zymd} and $\nu(\cdot)=Leb(\cdot\cap E)$.

(ii) For any step size $\eta\in(0,1)$, the Markov transition kernel $Q_\eta$ defined in \eqref{zyh} is irreducible and strongly aperiodic.
\end{proposition}

\begin{proof}
(i). Suppose that $E$ is a compact subset of $\R^d$ such that $Leb(E)>0$. Then for all $x\in E$ and $A\in\mathcal B(\R^d)$,
\begin{align}\label{436}
Q_\eta(x,A)=\int_Aq_\eta(x,\tl x)\dif \tl x\geq \int_{A\cap E}q_\eta(x,\tl x)\dif y\geq \inf_{x\in E,\tl x\in E}q_\eta(x,\tl x) Leb(A\cap E).
\end{align}
It is easy to see that the transition density function $q_{\eta}(x,\tl x)$ in \eqref{zymd} is a bivariate continuous function. Since  $E$ is a compact set, there exists a point $(x_0,y_0)$ on the boundary of $E$ such that
 \begin{align}\label{437}
\inf_{x\in E,\tl x\in E}q_\eta(x,\tl x)=q_\eta(x_0,y_0)\in(0,1].
\end{align}
Combining \eqref{436} and \eqref{437}, we have for all $x\in E$ and $A\in\mathcal B(\R^d)$
\begin{align*}
Q_\eta(x,A)\geq\epsilon \nu(A),
\end{align*}
where $\nu(\cdot)=Leb(\cdot\cap E)$ and
\begin{align*}
\epsilon=\inf_{x\in E,\tl x\in E}q_\eta(x,\tl x)\in(0,1].
\end{align*}
Particularly, $Q_\eta(x,E)\geq\epsilon Leb (E)>0$ for all $x\in\R^d$.
By Definition \ref{small} , $E$ is an accessible $(1,\epsilon\nu)$-small set.

(ii). It is known from the proof of  part (i) that any compact subset $E$ of $\R^d$ satisfying $Leb(E)>0$ is an accessible $(1,\mu)$-small set with $\mu(E)=\epsilon\nu(E)>0$. Then according to Definition \ref{irre} and Definition \ref{ape} (iii), $Q_\eta$ is irreducible and strongly aperiodic.
\end{proof}

For the Markov kernel $Q_\eta$, the following lemma shows that the Lyapunov function condition holds. The proof is deferred to  Appendix \ref{a0}.

\begin{lemma}\label{lya}
Let $V(x)=1+|x|^2$.  Under Assumption \ref{ass}, the Markov kernel $Q_\eta$ satisfies for any $x\in\R^d$
\begin{align*}
Q_\eta V(x)\leq \lambda(\eta) V(x)+b(\eta) 1_{D_\eta}(x),
\end{align*}
where $\lambda(\eta)=1-m\eta+(1+3L^2)\eta^2+4(2-1/K)L^4\eta^4
$, $b(\eta)=[m+2d+L^2(1-1/K)]\eta+4dL^2(2-1/K)\eta^3+4(2-1/K)\eta^4$ and $D_{\eta}=\left\{x:|x|\leq\sqrt{\frac{b(\eta)}{m\eta}}\right\}$.
\end{lemma}

\subsection{Proof of Main results}
In this subsection, we will prove Theorem \ref{ms} and Proposition \ref{ms2}. First we state two propositions which will be used in the proof of Theorem \ref{ms}.

 Given a Markov chain $(X_n)_{n\geq0}$ and its Markov kernel $P$ on $X\times\mathcal X$. For $A\in\mathcal X$, the first hitting time $\tau_A$ and return time $\sigma_A$ of the set $A$ by the Markov chain $(X_n)_{n\geq0}$ are defined respectively by
 \begin{align*}
 \tau_A&=\inf\{n\geq0:X_n\in A\},\\
 \sigma_A&=\inf\{n\geq1:X_n\in A\}.
 \end{align*}

\begin{proposition}[Proposition 4.3.3 in \cite{MC}]\label{sp}
Let $P$ be a Markov kernel on $X\times\mathcal X$ and $A\in\mathcal X$. If there exists a function $V:X\to[1,\infty]$, $\lambda\in[0,1)$ and $b< \infty$ such that $PV\leq\lambda V+b1_A$, then for all $x\in X$,
\begin{align*}
\E_x[\lambda^{-\sigma_A}]\leq V(x)+b\lambda^{-1},
\end{align*}
\end{proposition}

\begin{proposition}[Theorem 11.4.2 in \cite{MC}]\label{sp2}
Let $P$ be a Markov kernel on $X\times\mathcal X$. Assume that there exists an accessible $(1,\mu)$-small set $A$ and $\beta>1$ such that $\mu(A)>0$ and $\sup_{x\in A}\E_x[\beta^{\sigma_A}]<\infty$.  Then $P$ has a unique invariant probability measure $\pi_P$, and there exist $\delta>1$ and $\theta<\infty$ such that for all probability measure $\zeta$ on $\mathcal X$,
\begin{align*}
\sum_{k=1}^\infty \delta^kd_{\rm TV}(\zeta P^k, \pi_P)\leq \theta \E_{\zeta}[\beta^{\sigma_A}].
\end{align*}
\end{proposition}

With the above auxiliary results in hand, we are in a position to prove our main result.

\begin{proof}[Proof of Theorem \ref{ms}]

We have showed in Lemma \ref{lya} that the Markov kernel of $(\tl X_n)_{n\geq0}$ satisfies
\begin{align*}
Q_\eta V(x)\leq \lambda(\eta) V(x)+b(\eta) 1_{D_\eta}(x).
\end{align*}
It is obvious that there exists an $\eta_0$  such that when $\eta\in(0,\eta_0)$, $0<\lambda(\eta)<1$.

Since for $\eta\in(0,\eta_0)$, $D_\eta$ is a compact subset of $\R^d$ and ${\rm Leb}(D_\eta)>0$, from Proposition \ref{sms}, we have for any $\eta\in(0,\eta_0)$, $D_\eta$ is an accessible $(1,\mu)$-small set, and 
\begin{align*}
\mu(D_\eta)=\epsilon {\rm Leb}(D_\eta)>0.
\end{align*}
Then by Proposition \ref{sp}, for all $x\in\R^d$,
\begin{align*}
\E_x\left[\left(\frac{1}{\lambda(\eta)}\right)^{\sigma_{D_\eta}}\right]\leq V(x)+\frac{b(\eta)}{\lambda(\eta)}.
\end{align*}
So that
\begin{align*}
\sup_{x\in D_\eta}\E_x\left[\left(\frac{1}{\lambda(\eta)}\right)^{\sigma_{D_\eta}}\right]\leq \sup_{x\in D_\eta}\left(V(x)+\frac{b(\eta)}{\lambda(\eta)}\right)<\infty.
\end{align*}
According to Proposition \ref{sp2}, $(\tl X_n)_{n\geq0}$ has a  unique stationary distribution $\pi_\eta$, and there exist $\delta>1$ and $\theta<\infty$ such that for all initial state $X_0=x\in\R^d$,
\begin{align*}
\sum_{n=1}^\infty\delta^nd_{\rm TV}(\mathcal L(\tl X_n),\pi_\eta)
\leq&\theta\E_x\left[\left(\frac{1}{\lambda(\eta)}\right)^{\sigma_{D_\eta}}\right]\\
\leq&\theta\left(1+|x|^2+\frac{b(\eta)}{\lambda(\eta)}\right)\\
<&+\infty,
\end{align*}
which implies that 
\begin{align*}
d_{\rm TV}(\mathcal L(\tl X_n),\pi_\eta)=o(\delta^{-n}), \quad n\to\infty.
\end{align*}

By Lemma \ref{lya} again, 
\begin{align*}
Q_\eta (|x|^2)&\leq\left[1-2m\eta+(1+3L^2)\eta^2+4(2-1/K)L^4\eta^4\right]|x|^2\nonumber\\
&\qquad+\eta d\{2+4\eta^2L^2(2-1/K)+L^2(1-1/K)\}\nonumber
\end{align*}
When $\eta$ small enough, $1-2m\eta+(1+3L^2)\eta^2+4(2-1/K)L^4\eta^4\leq 1-m\eta<1$.

Inductively, for any $n\geq2$,
\begin{align}\label{10111201}
&Q_\eta^n(|x|^2)\nonumber \\
\leq& (1-m\eta)^n|x|^2+\eta d\{2+4\eta^2L^2(2-1/K)+L^2(1-1/K)\}\sum_{i=0}^{n-1}(1-m\eta)^i.
\end{align}
Let $n\to\infty$, the above inequality implies that 
\begin{align*}
\pi_\eta(|x|^2)\leq \frac{d}{m}\{2+4\eta^2L^2(2-1/K)+L^2(1-1/K)\}.
\end{align*}
The proof is complete.
\end{proof}

Next we turn to prove Proposition \ref{ms2}.

\begin{proof}[Proof of Proposition \ref{ms2}]
Recall the Langevin dynamics
\begin{align*}
\dif X_t=-\nabla U(X_t)\dif t+\sqrt 2\dif B_t
\end{align*}
and its Poisson randomized midpoint Euler-Maruyama scheme
\begin{align*}
\tl{X}_{1}=\tl{X}_{0}-\eta\nabla U(\tl{X}_{0})+\eta\sum_{i=0}^{K-1}H_i(\nabla U(\tl{X}_{0})-\nabla U(\hat{X}_i))+\sqrt{2}B_\eta,
\end{align*}
where
\begin{align*}
\hat{X}_i=\tl{X}_{0}-\frac{i\eta}{K}\nabla U(\tl{X}_{0})+\sqrt{\frac{2i\eta}{K}}\xi_i.
\end{align*}

Let $\tl{X}_{0}\sim\pi_\eta$ and $X_0\sim\pi$ be two independent random variables. Define 
\begin{align*}
X_t=X_0-\int_0^t\nabla U(X_s)\dif s+\sqrt 2B_t, \quad t\in[0,\eta]
\end{align*}
to be the solution of Langevin dynamics with initial value $X_0$. Therefore, $\tl{X}_{1}\sim\pi_\eta$ and $X_\eta\sim\pi$ are also independent and $\|\tl{X}_{1}-X_\eta\|_{2}=\|\tl{X}_{0}-X_0\|_{2}$.

Also, it easily follows
\begin{align*}
&X_\eta-\tl{X}_{1}\\
=&X_0-\int_0^\eta\nabla U(X_s)\dif s-\tl{X}_{0}+\eta \nabla U(\tl{X}_{0})-\eta\sum_{i=0}^{K-1}H_i(\nabla U(\tl{X}_{0})-\nabla U(\hat{X}_i))\\
=&X_0-\tl{X}_{0}-\eta\Big(\nabla U(X_0)-\nabla U(\tl{X}_{0})\Big)-\int_0^\eta\nabla U(X_s)-\nabla U(X_0)\dif s\\
&\quad-\eta\sum_{i=0}^{K-1}H_i(\nabla U(\tl{X}_{0})-\nabla U(\hat{X}_i))
\end{align*}
So we have
\begin{align}\label{tsg}
\|X_\eta-\tl{X}_{1}\|_{2}\leq
&\Big\|X_0-\tl{X}_{0}-\eta\Big(\nabla U(X_0)-\nabla U(\tl{X}_{0})\Big)\Big\|_{2}+\Big\|\int_0^\eta\nabla U(X_s)-\nabla U(X_0)\dif s\Big\|_{2}\nonumber\\
&+\eta\Big\|\sum_{i=0}^{K-1}H_i(\nabla U(\tl{X}_{0})-\nabla U(\hat{X}_i))\Big\|_{2}
\end{align}

Since $U$ is $L$-gradient Lipschitz and $m$-strongly convex and  suppose $\eta\in(0,\frac{2}{L+m})$,
\begin{align}\label{tsg1}
\Big\|X_0-\tl{X}_{0}-\eta\Big(\nabla U(X_0)-\nabla U(\tl{X}_{0})\Big)\Big\|_{2}\leq (1-m\eta)\|X_0-\tl{X}_{0}\|_{2}.
\end{align}

For the second term in \eqref{tsg}, note
\begin{align*}
\Big\|\int_0^\eta\nabla U(X_s)-\nabla U(X_0)\dif s\Big\|_{2}=&\sqrt{\E\left[\left|\int_0^\eta\nabla U(X_s)-\nabla U(X_0)\dif s\right|^2\right]}\\
\leq&\eta^{\frac{1}{2}}\sqrt{\int_0^\eta\E[|\nabla U(X_s)-\nabla U(X_0)|^2]\dif s}\\
\leq&\eta L\left\{\E\left[\sup_{0\leq t\leq \eta}\|X_t-X_0\|^2\right]\right\}^{\frac{1}{2}}.
\end{align*}
According to  Lemma 2 in \cite{hy},
\begin{align*}
\E\left[\sup_{0\leq t\leq \eta}\|X_t-X_0\|^2\right]\leq
4\eta^2\E[\|\nabla U(X_0)\|^2]+8L^2d\eta^3+2d\eta.
\end{align*}
Let $\mathcal L$ be the generator of Langevin dynamics \eqref{ld},
\begin{align*}
\mathcal L f=\Delta f-\langle \nabla U, \nabla f\rangle.
\end{align*}
Since $\nabla^2U\preceq L I_d$, then $\Delta  U\leq Ld$. Take $f=U$, we have
\begin{align*}
0=\E_\pi\mathcal L U=\E_\pi[\Delta  U-\|\nabla U\|^2],
\end{align*}
which implies
\begin{align}\label{1234}
\E[\|\nabla U(X_0)\|^2]\leq Ld.
\end{align}
 Therefore
\begin{align}\label{tsg2}
\Big\|\int_0^\eta\nabla U(X_s)-\nabla U(X_0)\dif s\Big\|_{2}\leq\eta L(2d\eta+4\eta^3Ld+8L^2\eta^3d)^{\frac{1}{2}}.
\end{align}

For the third term in \eqref{tsg}, note 
\begin{align*}
\E[|\nabla U(\tl{X}_{0})-\nabla U(\hat{X}_i)|^2]&\leq L^2\E[|\hat{X}_i-\tl{X}_{0}|^2]\\
&\leq L^2\E\left[\left|-\frac{i\eta}{K}\nabla U(\tl{X}_{0})+\sqrt{\frac{2i\eta}{K}}\xi_i\right|^2\right]\\
&\leq L^2\left[\frac{i^2\eta^2 L^2}{K^2}\pi_\eta(|x|^2)+\frac{2i\eta d}{K}\right].
\end{align*}

Then
\begin{align*}
&\E\left[\Big(\sum_{i=0}^{K-1}H_i(\nabla U(\tl{X}_{0})-\nabla U(\hat{X}_i))\Big)^2\right]\\
\leq &\E\left[2\left(\sum_{i=0}^{K-1}(H_i-\frac{1}{K})(\nabla U(\tl{X}_{0})-\nabla U(\hat{X}_i))\right)^2+2\left(\sum_{i=0}^{K-1}\frac{1}{K}(\nabla U(\tl{X}_{0})-\nabla U(\hat{X}_i))\right)^2\right]\\
\leq&2\frac{1}{K}(1-\frac{1}{K})\sum_{i=0}^{K-1}\E[|\nabla U(\tl{X}_{0})-\nabla U(\hat{X}_i)|^2]+2\frac{1}{K}\sum_{i=0}^{K-1}\E[|\nabla U(\tl{X}_{0})-\nabla U(\hat{X}_i)|^2]\\
\leq&(\frac{4}{K}-\frac{2}{K^2})L^2\sum_{i=0}^{K-1}\left[\frac{i^2\eta^2 L^2}{K^2}\pi_\eta(|x|^2)+\frac{2i\eta d}{K}\right]\\
\leq&(4-\frac{2}{K})L^2\left(\frac{1}{3}\eta^2 L^2\pi_\eta(|x|^2)+\eta d\right).
\end{align*}
 So we have 
\begin{align}\label{tsg3}
\eta\Big\|\sum_{i=0}^{K-1}H_i(\nabla U(\tl{X}_{0})-\nabla U(\hat{X}_i))\Big\|_{2}\leq \eta\sqrt{4-\frac{2}{K}} L\left(\frac{1}{3}\eta^2 L^2\pi_\eta(|x|^2)+\eta d\right)^{\frac{1}{2}}.
\end{align}

Combining \eqref{tsg}-\eqref{tsg3} together, we get
\begin{align*}
\|X_0-\tl{X}_{0}\|_{2}=\|X_\eta-\tl{X}_{1}\|_{2}\leq& (1-m\eta)\|X_0-\tl{X}_{0}\|_{2}+\eta L(2d\eta+4\eta^3Ld+8L^2\eta^3d)^{\frac{1}{2}}\\
&+ \eta\sqrt{4-\frac{2}{K}} L\left(\frac{1}{3}\eta^2L^2\pi_\eta(|x|^2)+\eta d\right)^{\frac{1}{2}},
\end{align*}
and then
\begin{align*}
\|X_0-\tl{X}_{0}\|_{2}\leq \frac{L}{m}\left[(2d\eta+4\eta^3Ld+8L^2\eta^3d)^{\frac{1}{2}}+\sqrt{4-\frac{2}{K}} \left(\frac{1}{3}\eta^2 L^2\pi_\eta(|x|^2)+\eta d\right)^{\frac{1}{2}}\right].
\end{align*}
The desired result now follows by the fact that $W_2(\pi_\eta,\pi)\leq \|X_0-\tl{X}_{0}\|_{2}$.

\end{proof}

\section{ Convergence rate  in $d_{\Gcal}$ distance for PRLMC}\label{2127}
In this section, we will use  the Lindeberg replacement 
  method developed in Pag\`es and Panloup \cite{pages} to prove Theorem \ref{ms3}. This method has been widely used in recent years to study the Euler-Maruyama approximations of ergodic stochastic differential equations, such as \cite{ch,  jtp2, de,  gu,ji}. We state several auxiliary lemmas and propositions  in Subsection \ref{sec:result}, whose rigorous  proofs will be provided in Appendix \ref{a1}. In Subsection \ref{1421}, we prove Theorem \ref{ms3}.

\subsection{Auxiliary results}\label{sec:result}

Let $X_t^{x}$  be the Langevin dynamics given by  \eqref{ld}. For any $0\leq s\leq t<\infty$, denote the semigroup of process $(X_t)_{t\geq 0}$ by $P_{s,t}$, that is, for any bounded Borel function $f\in \mathcal{B}_b(\R^d,\R)$ and any $x\in \R^d $, one has
\begin{eqnarray*}
P_{s,t}f(x) = \E\left[f(X_t)\Big|X_{s}=x\right].
\end{eqnarray*}
Notice that the process $(X_t)_{t\geq 0}$ is time homogeneous, we further denote
\begin{align*}
P_{t-s}f(x)=P_{s,t}f(x).
\end{align*}
For the discrete time non-homogeneous process $(Y_{t_k})_{k\in \mathbb{N}_0}$ in \eqref{bbc}, we denote
\begin{align*}
\tilde{P}_{t_{i},t_{j}}f(x)=\mathbb{E}\left[f(Y_{t_{j}})\Big|Y_{t_{i}}=x\right], \quad 0\leq i\leq j<\infty, f\in\mathcal{B}_b(\R^d\times \mathbb{S},\R)
\end{align*}
and the probability trasition
\begin{align*}
\tilde{P}_{\gamma_{k+1}} f(x)=\tilde{P}_{t_{k},t_{k+1}}f(x), \quad k\in\mathbb{N}_0.
\end{align*}
For the sake of narrative, let $(\tilde{Y}_{\gamma}^{x})$ denote the one step Poisson random midpoint Euler-Maruyama scheme with any step size $\gamma\in(0,\gamma_{1}]$ starting from $x$ at time 0, that is,
\begin{align}\label{onestep}
\tl Y_{\gamma}^x=x-\gamma\nabla U(x)+\gamma\sum_{i=0}^{K-1}H_{i}(\nabla U(x)-\nabla U(\hat Y_{i}))+\sqrt{2\gamma}\xi,
\end{align}
where 
\begin{align*}
\hat Y_{i}=x-\frac{i\gamma}{K}\nabla U(x)+\sqrt{\frac{2i\gamma}{K}}\xi_i,
\end{align*}
and then
\begin{eqnarray*}
\E f(\tilde{Y}_{\gamma}^{x})=\tilde{P}_{\gamma} f(x).
\end{eqnarray*}

First, we  give some useful moment estimates of   $(X_t)_{t\geq0}$ and  $(Y_{t_{k}})_{k\geq 0}$ .

\begin{lemma}\label{LYexp}
Under the above notations and  Assumptions \ref{ass}, for any initial state $x$ and $t\geq 0$, 
\begin{eqnarray}\label{Xmoment}
\mathbb{E}|X_{t}^{x}|^{2}\leq |x|^2e^{-2mt}+(d/m)(1-e^{-2mt}),
\end{eqnarray}
and for any $k\geq1$, there exists $C>0$ independent of $\Gamma=(\gamma_n)_{n\in \mathbb{N}}$ such that
\begin{eqnarray}\label{Yexp}
\mathbb{E}|Y_{t_{k}}^{x}|^{2}\leq C(1+|x|^{2}).
\end{eqnarray}
\end{lemma}

\begin{lemma}\label{one}
 Let $(\tilde{Y}_{\gamma}^{x,m})$ be the one step Poisson random midpoint Euler-Maruyama  scheme with step $\gamma\in(0,\gamma_{1}]$, $\gamma_1<1$. Then under Assumptions \ref{ass}, for every $t\in[0,\gamma]$ , there exists some positive constant $C$ independent of $t$ such that
\begin{align}
\|X_{t}^{x}-x\|_2&\leq C (1+|x|)t^{\frac{1}{2}},\label{Xdiff}\\
\|\tilde{Y}_{\gamma}^{x}-x\|_2&\leq C (1+|x|)\gamma^{\frac{1}{2}},\label{Ydiff}\\
\|X_{\gamma}^{x}-\tilde{Y}_{\gamma}^{x}\|_2&\leq C (1+|x|)\gamma^{\frac{3}{2}}\label{X-Ydiff}.
\end{align}
\end{lemma}

Next, we state a proposition about the regularity of Langevin  semigroup $(P_t)_{t\geq0}$.
\begin{proposition}\label{tidu}
Let $(P_t)_{t\geq0}$ be the semigroup of Langevin diffusion $(X_t)_{t\geq0}$. Under Assumption \ref{ass} , there exists some positive constant $C$ independent of $t$ such that for any $h\in\Gcal$ and $t\geq0$
\begin{align*}
\|\nabla P_th\| &\leq \| \nabla h \|e^{-mt},  \\
\|\nabla^2 P_th\|&\leq C (\| \nabla h \|+\| \nabla^2 h \|) e^{-\frac{m}{2}t}
\end{align*}
 where $m$ is in \eqref{m}.
\end{proposition}

\begin{proposition}\label{prop:onestep}
 Under Assumptions \ref{ass}, for any $x\in\R^d$, $0<\gamma\leq \gamma_1<1$ and $f\in\mathcal{C}^2(\R^d,\R)$, then there exists some positive constant $C$  independent of $\gamma$ such that 
\begin{align*}
|\E f(X_{\gamma}^{x})-\E f(\tilde{Y}_{\gamma}^{x})|\leq C(1+|x|^2)(\|\nabla f\|+\|\nabla f^{2}\|)\gamma^{2}.
\end{align*}
\end{proposition}

\begin{lemma}\label{Step}
Let $(\gamma_n)_{n\in \mathbb{N}}$ be a non-increasing positive sequence. If
\begin{align*}
\omega=\limsup_{n\rightarrow\infty} \{ (\gamma_{n}-\gamma_{n+1}) \gamma_{n+1}^{-2} \}<+\infty,
\end{align*}
let $m>2\omega$ and $(u_n)_{n\in \mathbb{N}}$ be a sequence defined by $u_0=0$ and
\begin{align*}
u_n=\sum_{k=1}^n \gamma_k^{2}e^{-\frac{m}{2}(t_n-t_k)}, \ \ \   n\geq1.
\end{align*}
Then, we have
\begin{align}\label{s1}
\limsup_{n\rightarrow\infty} \{ u_n \gamma_n^{-1} \}<+\infty,
\end{align}
and 
\begin{align}\label{exp1}
\lim_{n\rightarrow\infty} e^{-m t_{n}} \gamma_{n}^{-1}=0.
\end{align}
\end{lemma}
The proof of Lemma \ref{Step} follows the same lines as the proof of \cite[ Lemma 3.7]{jtp2} and is omitted.

\subsection{Proof of Theorem \ref{ms3}}\label{1421}

With the help of the above preliminary propositions and lemmas, we are in the position to prove Theorem \ref{ms3}.

\begin{proof}[Proof of Theorem \ref{ms3}]
For any $n\geq1$, by the definition of the  $d_\Gcal$, we have
\begin{align}\label{R1}
d_\Gcal \left(\Lcal\left( X_{t_n}^{x}\right),\Lcal\left( Y_{t_n}^{x}\right)\right) 
=&\sup_{h\in\Gcal}\{\E h(X_{t_n}^{x})-\E h(Y_{t_n}^{x}) \}\nonumber\\
=&\sup_{h\in\Gcal} \{P_{\gamma_1}\circ \cdots\circ P_{\gamma_n}h(x)-\tl{P}_{\gamma_1}\circ \cdots\circ\tl{P}_{\gamma_n}h(x) \} \nonumber\\
\leq&\sup_{h\in\Gcal}\left\{ \sum_{k=1}^n|\tl{P}_{\gamma_1}\circ \cdots\circ\tl{P}_{\gamma_{k-1}}\circ(P_{\gamma_{k}}-\tl{P}_{\gamma_k})\circ P_{t_n-t_k}h(x)| \right\}.
\end{align}
It follows from Propositions \ref{tidu} and \ref{prop:onestep} that for any $k \in\{1,2,\cdots,n\}$,
\begin{align*}
|(P_{\gamma_{k}}-\tl{P}_{\gamma_k})\circ P_{t_n-t_k}h(x)|
\leq& C(1+|x|^2)(\|\nabla P_{t_n-t_k}h\|+\|\nabla^2 P_{t_n-t_k}h\|)\gamma_k^{2}\\
\leq& C(1+|x|^2)(\|\nabla h\|+\|\nabla^2h\|)e^{-\frac{m}{2}(t_n-t_k)}\gamma_k^{2}.
\end{align*}
Then, integrating with respect to $\tl{P}_{\gamma_1}\circ \cdots\circ\tl{P}_{\gamma_{k-1}}$ yields
\begin{align}\label{R2}
&|\tl{P}_{\gamma_1}\circ \cdots\circ\tl{P}_{\gamma_{k-1}}\circ(P_{\gamma_{k}}-\tl{P}_{\gamma_k})\circ P_{t_n-t_k}h(x)|\nonumber \\
\leq&  C(\|\nabla h\|+\|\nabla^2h\|)e^{-\frac{m}{2}(t_n-t_k)}\gamma_k^{2}(1+\E|Y_{t_{k-1}}^{x}|^2)\nonumber\\
\leq&  C(1+|x|^2)(\|\nabla h\|+\|\nabla^2h\|)e^{-\frac{m}{2}(t_n-t_k)}\gamma_k^{2},
\end{align}
where the last inequality holds by (\ref{Yexp}). Hence, combining \eqref{R1} and \eqref{R2} together yields
\begin{align}\label{dh0}
d_\Gcal \left(\Lcal\left( X_{t_n}^{x}\right),\Lcal\left( Y_{t_n}^{x}\right)\right) &\leq  C(1+|x|^2)\sum_{k=1}^ne^{-\frac{m}{2}(t_n-t_k)}\gamma_k^{2}
\leq C(1+|x|^2)\gamma_n,
\end{align}
where the last inequality holds from inequality  \eqref{s1}.

Next we consider $d_\Gcal \left(\Lcal\left( X_{t_n}^{x}\right),\pi\right)$.  According to 
Durmus and Moulines \cite[Proposition 1 ]{du2}, under Assumption \ref{ass},
\begin{align*}
W_2\left(\Lcal\left( X_{t_n}^{x}\right),\pi\right)\leq e^{-mt_n}\{|x|+(d/m)^{1/2}\}.
\end{align*}

Denote
\begin{align*}
\mathcal F =\{f|f:\R^d\to\R,\|f\|_{\rm Lip}\leq 1\}.
\end{align*}
By the definition of the function class $\mathcal G$,  there exists a constant $M_1$   such that $\|\nabla h\|\leq M_1$ for any $h\in\mathcal G$. So we have
\begin{align}\label{dh1}
d_\Gcal \left(\Lcal\left( X_{t_n}^{x}\right),\pi\right)=\sup_{h\in\mathcal G}M_1\{\E \frac{h}{M_1}(X_{t_n}^{x})-\pi(\frac{h}{M_1})\}\leq M_1\sup_{h\in\mathcal F}\{\E h(X_{t_n}^{x})-\pi(h)\}.
\end{align}

 By Kantorovich dual formula \cite[Remark 6.5]{v}, we have
 \begin{align}\label{dh2}
&\sup_{h\in\mathcal F}\{\E h(X_{t_n}^{x})-\pi(h)\}\nonumber\\
=&W_1\left(\Lcal\left( X_{t_n}^{x}\right),\pi\right)\leq W_2\left(\Lcal\left( X_{t_n}^{x}\right),\pi\right)\leq e^{-mt_n}\{|x|+(d/m)^{1/2}\}.
\end{align}

Furthermore, combining \eqref{exp1}, \eqref{dh0}, \eqref{dh1} and  \eqref{dh2}  together yields 
\begin{align*}
 d_\Gcal \left(\Lcal\left( Y_{t_n}^{x}\right),\pi\right)
 \leq& d_\Gcal \left(\Lcal\left( X_{t_n}^{x}\right),\Lcal\left( Y_{t_n}^{x}\right)\right)+d_\Gcal \left(\Lcal\left( X_{t_n}^{x}\right),\pi\right)\\
 \leq&C\left(1+|x|^{2}\right)(\gamma_{n}+e^{-m t_{n}})\\
 \leq& C\left(1+|x|^{2}\right)\gamma_{n}.
\end{align*}
The proof is complete.
\end{proof}

\section{Non-asymptotic bounds in 2-Wasserstein distance for PRLMC}\label{sec5}

\subsection{Proofs of Theorem \ref{ms4}, Corollary \ref{ms5} and \ref{ms6}}


To show Theorem \ref{ms4}, it suffices to get some bounds on $W_2(\Lcal(Y_{t_n}),\pi P_{t_n})$ since $\pi$ is invariant for $P_t$ for all $t\geq0$. To do so,  we construct a coupling between the Langevin diffusion and the linear interpolation  of the PRLMC. Given a sequence of step size $\Gamma=(\gamma_n)_{n\geq1}$ , let $t_n=\sum_{i=1}^n\gamma_i$. 
For all $n\geq0$ and $t\in[t_n,t_{n+1})$, denote the synchronous coupling $(X_t, \bar Y_t)_{t\geq0}$ by
\begin{align}\label{1220tsg}
\left\{
\begin{aligned}
    &X_t=X_{t_n}-\int_{t_n}^t\nabla U(X_s)\dif s+\sqrt2(B_t-B_{t_n})\\
    &\bar Y_t=\bar Y_{t_n}-(t-t_n)\nabla U(\bar Y_{t_n})+(t-t_n)\sum\limits_{i=0}^{K-1}H_{n,i}(\nabla U(\bar Y_{t_n})-\nabla U(\hat Y_{t,i}))\\
    &\qquad\quad+\sqrt2(B_t-B_{t_n}),
\end{aligned}
\right.
\end{align}
where the midpoints in the interpolation process are  defined by
\begin{align}\label{1092103}
    \hat Y_{t,i}=\bar Y_{t_n}-\frac{i(t-t_n)}{K}\nabla U(\bar Y_{t_n})+\sqrt{\frac{2i}{K}}(B^\prime_{t}-B^\prime_{t_n}),
\end{align}
and $(B^\prime_{t})_{t\geq0}$ is a Brownian motion independent of $(B_{t})_{t\geq0}$. 
Let $(\mathcal F_t)_{t\geq0}$ be the filtration associated with $(X_0,Y_0)$, $(B_t)_{t\geq0}$ and $(B^\prime_t)_{t\geq0}$. 

Here are three stochastic processes: $(X_t)_{t\geq0}$ represents the Langevin dynamics, $(Y_{t_n})_{n\geq0}$
denotes the PRLMC, and  $(\bar Y_t)_{t\geq0}$ is the linear interpolation of PRLMC. Let $X_0\sim\pi$, $\bar Y_0=Y_0=x$,  we have that for all $n\geq1$, $Y_{t_n}\overset{d}{=}\bar Y_{t_n}$ and then $W_2^2(\Lcal(Y_{t_n}),\pi )=W_2^2(\Lcal(\bar Y_{t_n}),\pi P_{t_n} )\leq \E[|X_{t_n}-\bar Y_{t_n}|^2]$. We present the following lemma on the recurrence of $\E[|X_{t_n}-\bar Y_{t_n}|^2]$.

\begin{lemma}\label{s4}
Assume Assumptions \ref{ass} and \ref{ass2} hold. Let $(\gamma_k)_{k\geq1}$ be a non-increasing sequence with $\gamma_1\leq 1/(m+L)$. Let $(X_t,\bar Y_t)_{t\geq0}$  be defined by \eqref{1220tsg}, then almost surely for all $n\geq0$ and $\epsilon>0$,
\begin{align*}
    \E^{\mathcal F_{t_n}}[|X_{t_{n+1}}-\bar Y_{t_{n+1}}|^2]&\leq
    \{1-\gamma_{n+1}(\kappa-4\epsilon)\}|X_{t_{n}}-\bar Y_{t_{n}}|^2\\
    &\quad +\gamma_{n+1}^3\Big\{L^4[\gamma_{n+1}+(3\ep)^{-1}]|X_{t_{n}}|^2\\
    &\qquad\qquad\ +L^4[(8-4/K)\gamma_{n+1}+(3\ep)^{-1}]|\bar Y_{t_{n}}|^2\\
    &\qquad\qquad\ +L^2d[10-4/K+L^2\gamma_{n+1}^2/6+(4\ep)^{-1}L^2\gamma_{n+1}]\\
    &\qquad\qquad\ + 2(3\ep)^{-1}d^2\tl L^2\Big\},
\end{align*}
where $\kappa=2mL/(m+L)$.
\end{lemma}
\begin{proof}
    The proof is postponed to section \ref{pf51}.
\end{proof}

Now we turn to prove Theorem \ref{ms4}.
\begin{proof}[Proof of Theorem \ref{ms4}]
    Let $x\in\R^d$, $n\geq1$ and $\zeta_0=\pi\otimes\delta_x$. Let $(X_t,\bar Y_t)_{t\geq0}$ be the coupling processes defined by \eqref{1220tsg} with $(X_0,\bar Y_0)$ distributed according to $\zeta_0$.  By definition of $W_2$ and since for all $t\geq0$, $\pi$ is invariant for $P_t$, $W_2^2(\Lcal( Y_{t_n}),\pi)\leq \E[|X_{t_n}-\bar Y_{t_n}|^2]$. Lemma \ref{s4} with $\ep=\kappa/8$, Lemma \ref{LYexp} and \cite[Proposition 1]{du3} imply, using a straightforward induction, that for all 
    $n\geq0$,
\begin{align*}
    \E[|X_{t_n}-\bar Y_{t_n}|^2]\leq& u_n(\Gamma)\int_{\R^d}|y-x|^2\pi(\dif x)+A_n(\Gamma)\\
    \leq&2u_n(\Gamma)(|x|^2+d/m)+A_n(\Gamma)
\end{align*}
where $u_n(\Gamma)=\prod\limits_{k=1}^{n}(1-\kappa\gamma_k/2)$ and
\begin{align*}
    A_n(\Gamma)=&\sum_{i=1}^{n}\Big\{\gamma_i^3L^2d[10-4/K+L^2\gamma_{i}^2/6+2\kappa^{-1}L^2\gamma_{i}]+6\kappa^{-1}\gamma_i^3d^2\tl L^2\\
    &\qquad+ \delta_i\gamma_i^3L^4[(8-4/K)\gamma_i+3\kappa^{-1}]\\
    &\qquad+\tl\delta_i\gamma_i^3L^4(\gamma_i+3\kappa^{-1})\Big\}\prod_{k=i+1}^n(1-\kappa\gamma_k/2)
\end{align*}
with $\delta_i=\E[|\bar Y_{t_{i-1}}|^2]=\E[|Y_{t_{i-1}}|^2]\leq C(1+|x|^2)$ and $\tl\delta_i=\E[|X_{t_{i-1}}|^2]\leq d/m$, which hold according to \eqref{Yexp} and \cite[Proposition 1]{du3} respectively. The proof is complete.
\end{proof}

\begin{proof}[Proof of Corollary \ref{ms5} ]
    The proof follows the same line as that of \cite[Corollary 6]{du2} and is omitted.
\end{proof}

\begin{proof}[Proof of Corollary \ref{ms6} ]
    First note that the factor $C(1+|x|^2)$ is the upper bound for $\E[|Y_{t_{i-1}}|^2]$ in the decreasing step size setting, when we set $\gamma_k=\eta$ for all $k\geq1$, it should be replaced by the bound in \eqref{10111201}: for $i\geq1$
    \begin{align*}
        \E[|Y_{i\eta}|^2]\leq&(1-m\eta)^i|x|^2+\eta d\{2+4\eta^2L^2(2-1/K)+L^2(1-1/K)\}\sum_{j=0}^{i-1}(1-m\eta)^j\\
        \leq&(1-m\eta)^i|x|^2+ d\{2+4\eta^2L^2(2-1/K)+L^2(1-1/K)\}/m.
    \end{align*}

Then by Theorem \ref{ms4}, we have
\begin{align}\label{10111230}
    W_2^2(\Lcal(Y_{n\eta}),\pi)\leq 2(1-\kappa\eta/2)^n\{|x|^2+d/m\}+u_n(\eta),
\end{align}
where
\begin{align*}
    u_n(\eta)=&\eta^3\sum_{i=1}^{n}\Big\{L^2d[10-4/K+L^2\eta^2/6+2\kappa^{-1}L^2\eta]+6\kappa^{-1}d^2\tl L^2\\
    &\qquad+\Big [(1-m\eta)^{i-1}|x|^2+ d\{2+4\eta^2L^2(2-1/K)+L^2(1-1/K)\}/m\Big]\\
    &\qquad \quad \cdot L^4[(8-4/K)\eta+3\kappa^{-1}]\\
    &\qquad+dL^4(\eta+3\kappa^{-1})/m\Big\}(1-\kappa\eta/2)^{n-i}.
\end{align*}

It is easy to check that for $\eta\in(0,1/(m+L)]$ and any $x$,
\begin{align*}
   \lim_{n\to\infty}\sum_{i=1}^{n}(1-m\eta)^{i-1}(1-\kappa\eta/2)^{n-i}|x|^2=0.
\end{align*}

By Theorem \ref{ms} and the definition of total variation distance, we have that
 for any measurable function $f:\R^d\to\R$ with $\|f\|\leq1$,
\begin{align*}
    |\E f(Y_{n\eta})-\pi_\eta(f)|\to0, \quad n\to\infty.
\end{align*}
Then for any $h\in \mathcal C_b(\R^d;\R)$, 
\begin{align}\label{1407w}
    \E h(Y_{n\eta})-\pi_\eta(h)\to0, \quad n\to\infty.
\end{align}
By \cite[Page 7]{bp}, \eqref{1407w} means $\Lcal(Y_{n\eta})$ converges weakly to $\pi_\eta$, then it follows from \cite[Remark 6.12]{v} that
\begin{align*}
    W_2^2(\pi_\eta,\pi)\leq \liminf _{n\to\infty}W_2^2(\Lcal(Y_{n\eta}),\pi).
\end{align*}
Lastly, letting $n\to\infty$ in \eqref{10111230} completes the proof.

\end{proof}

\subsection{Proof of Lemma \ref{s4}}\label{pf51}
If Assumption \ref{ass} holds, then \cite[Theorem 2.1.12, Theorem 2.1.19]{ne} show that  for all $x,y\in\R^d$:
\begin{align}\label{S1}
    \langle\nabla U(y)-\nabla U(x),y-x\rangle\geq
    \frac{\kappa}{2}|y-x|^2+\frac{1}{m+L}|\nabla U(y)-\nabla U(x)|^2,
\end{align}
where 
\begin{align*}
    \kappa=\frac{2mL}{m+L}.
\end{align*}

\begin{proof}[Proof of Lemma \ref{s4}]
Set $\Theta_n=X_{t_n}-\bar Y_{t_n}$. By definition, we have 
\begin{align}\label{10091309}
    |\Theta_{n+1}|^2
\leq&|\Theta_{n}|^2+2\left|\int_{t_n}^{t_{n+1}}\nabla U(X_s)-\nabla U(\bar Y_{t_n})\dif s\right|^2\nonumber\\
    &+2\gamma_{n+1}^2\Big|\sum\limits_{i=0}^{K-1}H_{n,i}(\nabla U(\bar Y_{t_n})-\nabla U(\hat Y_{t_{n+1},i}))\Big|^2\nonumber\\
    &-2\langle \Theta_{n}, \int_{t_n}^{t_{n+1}}\nabla U(X_s)-\nabla U(\bar Y_{t_n})\dif s\rangle\nonumber\\
    &+2\langle \Theta_{n}, \gamma_{n+1}\sum\limits_{i=0}^{K-1}H_{n,i}(\nabla U(\bar Y_{t_n})-\nabla U(\hat Y_{t_{n+1},i}))\rangle.
\end{align}

For the second term, 
\begin{align}\label{10091400}
&2\left|\int_{t_n}^{t_{n+1}}\nabla U(X_s)-\nabla U(\bar Y_{t_n})\dif s\right|^2\nonumber\\
\leq&2\gamma_{n+1}^2|\nabla U(X_{t_n})-\nabla U(\bar Y_{t_n})|^2+2\gamma_{n+1}\int_{t_n}^{t_{n+1}}|\nabla U(X_s)-\nabla U(X_{t_n})|^2\dif s.
\end{align}

For the fourth term, it follows from \eqref{S1} that
\begin{align}\label{10091405}
&-2\langle \Theta_{n}, \int_{t_n}^{t_{n+1}}\nabla U(X_s)-\nabla U(\bar Y_{t_n})\dif s\rangle\nonumber\\
=&-2\gamma_{n+1}\langle \Theta_{n}, \nabla U(X_{t_n})-\nabla U(\bar Y_{t_n})\rangle-2\int_{t_n}^{t_{n+1}}\langle \Theta_{n}, \nabla U(X_s)-\nabla U(X_{t_n})\rangle\dif s\nonumber\\
\leq&-2\gamma_{n+1}\Big(\frac{\kappa}{2}|\Theta_{n}|^2+\frac{1}{m+L}|\nabla U(X_{t_n})-\nabla U(\bar Y_{t_n})|^2\Big)\nonumber\\
&\qquad\qquad\qquad\qquad\qquad\qquad\qquad\quad-2\int_{t_n}^{t_{n+1}}\langle \Theta_{n}, \nabla U(X_s)-\nabla U(X_{t_n})\rangle\dif s.
\end{align}

Plugging \eqref{10091400}, \eqref{10091405} into \eqref{10091309} and using the condition $\gamma_1\leq 1/(m+L)$, we obtain
\begin{align}\label{10091454}
    \E^{\mathcal F_n}[ |\Theta_{n+1}|^2]\leq &(1-\kappa\gamma_{n+1})|\Theta_{n}|^2+\underbrace{2\gamma_{n+1}\int_{t_n}^{t_{n+1}} \E^{\mathcal F_n}[|\nabla U(X_s)-\nabla U(X_{t_n})|^2]\dif s}_{=:A_1}\nonumber\\
    &+\underbrace{2\gamma_{n+1}^2 \E^{\mathcal F_n}\Big[\Big|\sum\limits_{i=0}^{K-1}H_{n,i}(\nabla U(\bar Y_{t_n})-\nabla U(\hat Y_{t_{n+1},i}))\Big|^2\Big]}_{=:A_2}\nonumber\\   &+\underbrace{2\int_{t_n}^{t_{n+1}}|\langle\Theta_n, \E^{\mathcal F_n}[\nabla U(X_s)-\nabla U(X_{t_n})]\rangle|\dif s}_{=:A_3}\nonumber\\ 
&+\underbrace{2\gamma_{n+1}\Big|\langle\Theta_{n},\frac{1}{K}\sum_{i=0}^{K-1}\E^{\mathcal F_n}[\nabla U(\bar Y_{t_n})-\nabla U(\hat Y_{t_{n+1},i})]\rangle\Big|}_{=:A_4}.
\end{align}

\textit{Term $A_1$}. By \cite[Lemma S2]{du3},
\begin{align}
    A_1\leq 2\gamma_{n+1}L^2\Big(d\gamma_{n+1}^2+dL^2\gamma_{n+1}^4/12+(1/2)L^2\gamma_{n+1}^3|X_{t_n}|^2\Big)
\end{align}

\textit{Term $A_2$}. By \eqref{tsg3}, 
\begin{align}
    A_2\leq \gamma_{n+1}^2(8-4/K)L^2(\gamma_{n+1}^2L^2|\bar Y_{t_n}|^2+\gamma_{n+1}d).
\end{align}

\textit{Term $A_3$}. Using It\^{o}'s formula, we have for all $s\in[t_n,t_{n+1})$,
\begin{align*}
\nabla U(X_s)-\nabla U(X_{t_n})=&\int_{t_n}^s\{-\nabla^2U(X_r)\nabla U(X_r)+\vec{\Delta}(\nabla U)(X_r)\}\dif r\nonumber\\
&+\sqrt2\int_{t_n}^s\nabla^2U(X_r)\dif B_r.
\end{align*}
By \cite[Eq. (S21), (S22)]{du3} and  for all $\epsilon>0$, $2|\langle a,b\rangle|\leq2\epsilon|a|^2+(2\epsilon)^{-1}|b|^2$, 
we have
\begin{align}
    A_3\leq &2\epsilon\gamma_{n+1}|\Theta_{n}|^2\nonumber\\
    &+(2\epsilon)^{-1}\int_{t_n}^{t_{n+1}}\Big|\E^{\mathcal F_n}\int_{t_n}^s\{-\nabla^2U(X_r)\nabla U(X_r)+\vec{\Delta}(\nabla U)(X_r)\}\dif r\Big|^2\dif s\nonumber\\
    \leq&2\epsilon\gamma_{n+1}|\Theta_{n}|^2+(2\epsilon)^{-1}\Big(2L^4\gamma_{n+1}^3|X_{t_n}|^2/3+L^4d\gamma_{n+1}^4/2+2\gamma_{n+1}^3d^2\tl L^2/3\Big)
\end{align}

\textit{Term $A_4$}. Consider $\nabla U(\hat Y_{t_{n+1},i})-\nabla U(\bar Y_{t_n})$, where $ \hat Y_{t_{n+1},i}$ is defined by \eqref{1092103} with $t=t_{n+1}$. It follows from It\^{o}'s formula that
\begin{align*}
    &\E^{\mathcal F_n}[\nabla U(\hat Y_{t_{n+1},i})-\nabla U(\bar Y_{t_n})]\\
    =&\int_{t_n}^{t_{n+1}}\E^{\mathcal F_n}\{-\frac{i}{K}\nabla^2U(\hat Y_{t,i})\nabla U(\bar Y_{t_n})+\frac{i}{K}\vec{\Delta}(\nabla U)(\hat Y_{t,i})\}\dif t,
\end{align*}
which together with the Cauchy-Schwarz inequality implies that
\begin{align*}
    &|\E^{\mathcal F_n}[\nabla U(\hat Y_{t_{n+1},i})-\nabla U(\bar Y_{t_n})]|^2\\
\leq&\gamma_{n+1}\int_{t_n}^{t_{n+1}}\Big|\E^{\mathcal F_n}\{-\frac{i}{K}\nabla^2U(\hat Y_{t,i})\nabla U(\bar Y_{t_n})+\frac{i}{K}\vec{\Delta}(\nabla U)(\hat Y_{t,i})\}\Big|^2\dif t\\
    \leq&2\gamma_{n+1}\frac{i^2}{K^2}\int_{t_n}^{t_{n+1}}\{L^4|\bar Y_{t_n}|^2+d^2\tl L^2\}\dif t\\
    =&2\gamma_{n+1}^2\frac{i^2}{K^2}(L^4|\bar Y_{t_n}|^2+d^2\tl L^2).
\end{align*}
Then 
\begin{align}\label{A4}
A_4\leq&2\epsilon\gamma_{n+1}|\Theta_n|^2+\gamma_{n+1}(2\epsilon)^{-1}\frac{1}{K}\sum_{i=0}^{K-1}|\E^{\mathcal F_n}[\nabla U(\hat Y_{t_{n+1},i})-\nabla U(\bar Y_{t_n})]|^2\nonumber\\
\leq&2\epsilon\gamma_{n+1}|\Theta_n|^2+\epsilon^{-1}\gamma_{n+1}^3(L^4|\bar Y_{t_n}|^2+d^2\tl L^2)/3.
\end{align}

Combining \eqref{10091454}-\eqref{A4} completes the proof.

\end{proof}

\begin{appendix}
\section{ Proof of the Lyapunov condition  }\label{a0}

\begin{proof}[Proof of Lemma \ref{lya}]
Recall the definition of the Markov chain $(\tl X_{k})_{k\geq0}$ in \eqref{prlmc}, for any initial state $\tl X_0=x\in\R^d$, 
\begin{align}\label{1955}
Q_\eta V(x)&=\E_x[V(\tl X_1)]\nonumber\\
&=\E[|\tl X_1|^2]+1\nonumber\\
&\leq1+|x|^2+2\eta^2|\nabla U(x)|^2+2\eta^2\E\left[\Big(\sum\limits_{i=0}^{K-1}H_i(\nabla U(x)-\nabla U(\hat X_i))\Big)^2\right]+2\eta d\nonumber\\
&\quad +2\langle x, -\eta\nabla U(x)\rangle+2\E\left[\left\langle x,\eta\sum\limits_{i=0}^{K-1}H_i(\nabla U(x)-\nabla U(\hat X_i))\right\rangle\right]\nonumber\\
&\leq1+|x|^2+2\eta^2L^2|x|^2+2\eta d+2\eta^2\E\left[\Big(\sum\limits_{i=0}^{K-1}H_i(\nabla U(x)-\nabla U(\hat X_i))\Big)^2\right]\nonumber\\
&\quad +2\E\left[\left\langle x,-\eta\sum\limits_{i=0}^{K-1}H_i \nabla U(\hat X_i))\right\rangle\right],
\end{align}
where $
\hat{X}_i=x-\frac{i\eta}{K}\nabla U(x)+\sqrt{\frac{2i\eta}{K}}\xi_i$, $i=0,1,\ldots,K-1$.

By \eqref{tsg3}, we have
\begin{align}\label{10101556}
    2\eta^2\E\left[\Big(\sum\limits_{i=0}^{K-1}H_i(\nabla U(x)-\nabla U(\hat X_i))\Big)^2\right]\leq 4\eta^2L^2(2-1/K)(\eta^2L^2|x|^2/3+\eta d).
\end{align}

For the last term in \eqref{1955}, note
\begin{align*}
2\E\left[\Big\langle x,-\eta\sum\limits_{i=0}^{K-1}H_i\nabla U(\hat X_i)\Big\rangle\right]=\frac{1}{K}
\sum\limits_{i=0}^{K-1}\langle -2\eta x, \E \nabla U(\hat X_i)\rangle.
\end{align*}
By Taylor's formula, there exists a random point $y_i$ depending on $x$ and $\hat X_i$ such that 
\begin{align*}
\nabla U(\hat X_i)&=\nabla U(x)+\nabla^2U(y_i)(\hat X_i-x)\\
&=\nabla U(x)+\nabla^2U(y_i)\Big(-\frac{i\eta}{K}\nabla U(x)+\sqrt{\frac{2i\eta}{K}}\xi_i\Big).
\end{align*}
Then
\begin{align*}
\langle -2\eta x,\E\nabla U(\hat X_i)\rangle=&\langle  -2\eta x, \nabla U(x)\rangle +\E\left[\Big\langle  -2\eta x, \nabla ^2 U(y_i)\Big(-\frac{i\eta}{K}\nabla U(x)\Big)\Big\rangle\right]\\
&\quad +\E\left[\Big\langle  -2\eta x, \nabla ^2 U(y_i)\sqrt{\frac{2i\eta}{K}}\xi_i\Big\rangle\right]\\
\leq&-2m\eta|x|^2+\frac{2i\eta^2L^2}{K}|x|^2+\eta^2|x|^2+\frac{2L^2i\eta d}{K}
\end{align*}
So we have
\begin{align}\label{bw3}
2\E\left[\Big\langle x,-\eta\sum\limits_{i=0}^{K-1}H_i\nabla U(\hat X_i)\Big\rangle\right]\leq (-2m\eta+L^2\eta^2+\eta^2)|x|^2+\frac{K-1}{K}L^2d\eta.
\end{align}

Combining \eqref{1955}-\eqref{bw3}, we obtain
\begin{align}\label{2217}
Q_\eta V(x)&\leq\left[1-2m\eta+(1+3L^2)\eta^2+4(2-1/K)L^4\eta^4\right]|x|^2\nonumber\\
&\qquad+1+\eta d\{2+4\eta^2L^2(2-1/K)+L^2(1-1/K)\}\nonumber\\
&\leq \lambda(\eta) V(x)+b(\eta)1_{D_\eta}(x)
\end{align}
where $\lambda(\eta)=1-m\eta+(1+3L^2)\eta^2+4(2-1/K)L^4\eta^4
$, $b(\eta)=[m+2d+L^2(1-1/K)]\eta+4dL^2(2-1/K)\eta^3+4(2-1/K)\eta^4$ and $D_{\eta}=\left\{x:|x|\leq\sqrt{\frac{b(\eta)}{m\eta}}\right\}$.
The proof is complete.
\end{proof}

\section{ Proofs of the auxiliary results in Section \ref{2127} and \ref{sec5}}\label{a1}

\begin{proof}[Proof of Lemma \ref{LYexp}]
(i) Recall the generator $\mathcal A$ of Langevin diffusion is
\begin{align*}
\mathcal L g=\Delta g-\langle \nabla U, \nabla g\rangle.
\end{align*}
Let $g(x)=|x|^2$, by \eqref{u2}, we have
\begin{align*}
\mathcal L g(x)=&2d-\langle \nabla U(x),2x\rangle\\
\leq&-2m|x|^2+2d.
\end{align*}
It follows from \cite[proof of Lemma 7.2]{gurvich2014diffusion} that
\begin{align}\label{e:Vm}
\mathbb{E} |X_t^x|^2\leq e^{-2mt }|x|^2+\frac{d(1-e^{-2mt})}{m},
\quad \forall t\geq 0.
\end{align}

(ii) By \eqref{2217}, there exist positive constants $C_{1},C_{2}$ such that
\begin{align*}
\mathbb{E}|Y_{t_{k+1}}|^{2}\leq&\left(1-2m\gamma_{k+1}+C_{1}\gamma_{k+1}^{2}\right)\mathbb{E}|Y_{t_{k}}|^{2}+C_{2}\gamma_{k+1}.
\end{align*}
Since $\lim_{k\rightarrow\infty}\gamma_{k}=0$, there exists $k_{0}\in\mathbb{N}$ such that for every  $k\geq k_{0}$ satisfying $C_{1}\gamma_{k+1}^{2} \leq m\gamma_{k+1}$ and
 \begin{align*}
\mathbb{E}|Y_{t_{k+1}}|^{2}\leq&\left(1-m\gamma_{k+1}\right)\mathbb{E}|Y_{t_{k}}|^{2}+C_{2}\gamma_{k+1},
\end{align*}
and $1-m\gamma_{k+1}>0$. Inductively, for any $k>k_{0}$, we have
\begin{align*}
\mathbb{E}|Y_{t_{k}}|^{2}\leq\prod_{_{j=k_{0}+1}}^{k}\left(1-m\gamma_{j}\right)\mathbb{E}|Y_{t_{k_{0}}}|^{2}
+C_{2}\sum_{j=k_{0}+1}^{k}\left[\gamma_{j}\prod_{l=j+1}^{k}\left(1-m\gamma_{l}\right)\right],
\end{align*}
where we used the convention that $\prod_{l=k+1}^{k}\left(1-m\gamma_{l}\right)=1$.
Note that
\begin{align*}
m\sum_{j=k_{0}+1}^{k}\gamma_{j}\prod_{l=j+1}^{k}\left(1-m\gamma_{l}\right)=1-\prod_{l=k_{0}+1}^{k}\left(1-m\gamma_{l}\right)\leq1,
\end{align*}
 Hence, the above two inequalities imply
\begin{align}\label{induc1}
\mathbb{E}|Y_{t_{k}}|^{2}\leq\mathbb{E}|Y_{t_{k_{0}}}|^{2}
+\frac{C_{2}}{m}, \quad \forall k\geq k_0.
\end{align}
In addition, by a standard argument (see Lamberton and Pag\`es \cite[Lemma 2]{LP02}), it is easy to verify that for every $0\leq k\leq k_{0}$,
\begin{align}\label{induc2}
\mathbb{E}|Y_{t_{k}}^{x}|^{2}\leq C(1+|x|^{2}).
\end{align}
Combining (\ref{induc1}) and (\ref{induc2}), for any $k\geq0$, we have
\begin{align*}
\mathbb{E}|Y_{t_{k}}^{x}|^{2}\leq C(1+|x|^{2}).
\end{align*}
The proof is complete.
\end{proof}

\begin{proof}[Proof of Lemma \ref{one}]
(i) First recall that
\begin{align*}
X_t^x=x-\int_0^t\nabla U(X_s)\dif s+\sqrt2B_t,
\end{align*}
it follows from Minkowski's inequality, \eqref{u1} and \eqref{Xmoment} that 
\begin{align*}
\|X_t^x-x\|_2=&\left\|\int_0^t\nabla U(X_s)\dif s\right\|_2+\sqrt{2t}\|B_1\|_2\\
\leq&\sqrt t\left\{\int_0^t\E[|\nabla U(X_s)|^2]\dif s\right\}^{\frac{1}{2}}+\sqrt{2td}\\
\leq&C(1+|x|)t+\sqrt{2td}\\
\leq&C(1+|x|)t^{\frac{1}{2}}.
\end{align*}

(ii) By the definition of $\tl Y_{\gamma}^x$ in \eqref{onestep},
\begin{align*}
\|\tl Y_{\gamma}^x-x\|_2\leq& \gamma\left\|\Big(\sum_{i=0}^{K-1}H_{i}-1\Big)\nabla U(x)\right\|_2+\gamma\left\|\sum_{i=0}^{K-1}H_{i}\nabla U(\hat Y_{i}))\right\|_2+\sqrt{2\gamma d}\\
\leq&\gamma \sqrt{1-\frac{1}{K}}L|x|+\sqrt{2\gamma d}+\gamma\left\|\sum_{i=0}^{K-1}H_{i}\nabla U(\hat Y_{i}))\right\|_2.
\end{align*}
By the Cauchy-Schwarz inequality  we have 
\begin{align*}
\E\left[\Big(\sum_{i=0}^{K-1}H_{i}\nabla U(\hat Y_{i}))\Big)^2\right]\leq&K\sum_{i=0}^{K-1}\E[H_i^2|\nabla U(\hat Y_{i}))|^2]\\
=&\sum_{i=0}^{K-1}\E[|\nabla U(\hat Y_{i}))|^2]\\
\leq&\frac{1}{2}KL^2|x|^2+\frac{1}{2}\gamma L^2dK.
\end{align*}

Combining the above inequalities, we obtain 
\begin{align*}
\|\tl Y_{\gamma}^x-x\|_2\leq C(1+|x|)\gamma^{\frac{1}{2}}.
\end{align*}

(iii) By Minkowski's inequality, \eqref{tsg2} and \eqref{tsg3}, there exists some positive constant $C$ such that 
\begin{align*}
\|X_{\gamma}^{x}-\tilde{Y}_{\gamma}^{x}\|_2&\leq
\Big\|\int_0^\gamma\nabla U(X_s)-\nabla U(x)\dif s\Big\|_{2}+\gamma\Big\|\sum_{i=0}^{K-1}H_i(\nabla U(x)-\nabla U(\hat{Y}_i))\Big\|_{2}\\
&\leq C(1+|x|)\gamma^{\frac{3}{2}}
\end{align*}
\end{proof}

To  prove Proposition \ref{tidu}, we  consider the derivative of $X_t^x$ with respect to the initial value $x$, which is called the Jacobi flow. Let $v\in\R^d$, the Jacobi flow $\nabla_vX_t^x$ along the direction $v$ is defined by
	\begin{align}\label{jf1}
		\nabla_vX_t^x=\lim_{\e\to 0}\frac{X_t^{x+\e v}-X_t^x}{\e}, \ t\geq0.
	\end{align}
	
	Similarly, for $v_1, v_2\in\R^d$, we can define
	
	\begin{align*}
		\nabla_{v_2}\nabla_{v_1}X_t^x=\lim_{\e\to 0}\frac{\nabla_{v_1}X_t^{x+\e v_2}-\nabla_{v_1}X_t^x}{\e}, \ t\geq0,
	\end{align*}

	We need the following lemma on the moment estimates of Jacobi flows.
	
	\begin{lemma}\label{jm}
		Under Assumption \ref{ass} and \ref{ass2}, for any $v_1, v_2, x\in\R^d$ and $t\geq0$, there exist a positive  constant $C_1$ depending on $m$ and $\|\nabla^3U\|_{\mathrm{op}}$ such that
		\begin{align*}
			\E[|\nabla_{v_1}X_t^x|^{4}]\leq &|v_1|^{4}e^{-4m t},\\
			\E[|\nabla_{v_2}\nabla_{v_1}X_t^x|^{2}]\leq &C_1 |v_1|^2|v_2|^2e^{-m t}.
		\end{align*}
		
	\end{lemma}

\begin{proof}
(i) It follows from \eqref{ld} that
		\begin{align*}
			X_t^x=x+\int_0^t -\nabla U(X_s^x)\dif s+\sqrt2 B_t.
		\end{align*}
		Combining with \eqref{jf1},  for any $v_1\in \R^d$, one has $\nabla_{v_1}X_0^x=v_1$ and
		\begin{align*}
			\dif\nabla_{v_1}X_t^x=-\nabla^2 U(X_t^x) \nabla_{v_1}X_t^x\dif t.
		\end{align*}
		
		Using It\^{o}'s formula for function $f(x)=|x|^{4}$ , one has
		\begin{align*}
			\dif|\nabla_{v_1}X_t^x|^{4}=&\langle 4|\nabla_{v_1}X_t^x|^{2}\nabla_{v_1}X_t^x,- \nabla^2U(X_t^x)\nabla_{v_1}X_t^x\rangle\dif t\\
			\leq&-4m|\nabla_{v_1}X_t^x|^{4}\dif t,
		\end{align*}
		which implies that
		\begin{align}\label{m1}
			\E[|\nabla_{v_1}X_t^x|^{4}]\leq |v_1|^{4}e^{-4m t}.
		\end{align}
		
		(ii) For any $v_1, v_2\in \R^d$, we know $\nabla_{v_2}\nabla_{v_1}X_0^x=0$ and
		\begin{align*}
			\dif \nabla_{v_2}\nabla_{v_1}X_t^x=\big(-\nabla ^3U(X_t^x)\nabla_{v_2}X_t^x\nabla_{v_1}X_t^x-\nabla^2U(X_t^x)\nabla_{v_2}\nabla_{v_1}X_t^x\big)\dif t
		\end{align*}
		
		Let $f(x)=|x|^2$,  using It\^{o}'s formula, we have
		\begin{align*}
			&\dif |\nabla_{v_2}\nabla_{v_1}X_t^x|^2\\
			=&\langle 2\nabla_{v_2}\nabla_{v_1}X_t^x,-\nabla^3U(X_t^x)\nabla_{v_2}X_t^x\nabla_{v_1}X_t^x-\nabla^2U(X_t^x)\nabla_{v_2}\nabla_{v_1}X_t^x\rangle\dif t\\
			\leq&\big(-2m|\nabla_{v_2}\nabla_{v_1}X_t^x|^{2}+2\|\nabla^3U\|_{\mathrm{op}}|\nabla_{v_2}\nabla_{v_1}X_t^x||\nabla_{v_1}X_t^x||\nabla_{v_2}X_t^x|\big)\dif t,
		\end{align*}
		
		By Young's inequality, for any positive constant $c$, one has
		\begin{align*}
			&c|\nabla_{v_2}\nabla_{v_1}X_t^x|\frac{1}{c}|\nabla_{v_1}X_t^x||\nabla_{v_2}X_t^x|
			\leq\frac{c^2|\nabla_{v_2}\nabla_{v_1}X_t^x|^{2}}{2}+\frac{|\nabla_{v_1}X_t^x|^2|\nabla_{v_2}X_t^x|^2}{2c^2}.
		\end{align*}
		
		Choose $c$ small enough such that $c^2\|\nabla^3U\|_{\mathrm{op}}<m$, we obtain
		\begin{align*}
			\dif |\nabla_{v_2}\nabla_{v_1}X_t^x|^2\leq \big(-m|\nabla_{v_2}\nabla_{v_1}X_t^x|^{2}+C|\nabla_{v_1}X_t^x|^2|\nabla_{v_2}X_t^x|^2\big)\dif t,
		\end{align*}
		where the constant $C$  depends on $m$ and $\|\nabla^3U\|_{\mathrm{op}}$.

		It follows from the comparison theorem (see \cite[Theorem 54]{Plotter2003Stochastic}) that
		\begin{align}\label{m2}
			\E[|\nabla_{v_2}\nabla_{v_1}X_t^x|^2]\leq &Ce^{-m t}\int_0^te^{m s}\E[|\nabla_{v_1}X_t^x|^2|\nabla_{v_2}X_t^x|^2]\dif s\nonumber\\
			\leq&C|v_1|^2|v_2|^2e^{-m t}\int_0^te^{-3m s}\dif s\nonumber\\
			\leq& C_1|v_1|^2|v_2|^2e^{-m t},
		\end{align}
		where we used the  moment estimate of one-order Jacobi flow in the second inequality and the constant $C_1$  depends on $m$ and $\|\nabla^3U\|_{\mathrm{op}}$.
\end{proof}

\begin{proof}[Proof of Proposition \ref{tidu}]
 For any $h\in\Gcal$ ,we have
\begin{align*}
\|\nabla P_th\|&=\sup_{x\in\R^d}\|\nabla P_th(x)\|_{\rm op}\\
&=\sup_{x\in\R^d,|u|=1}|\nabla_{u} P_th(x)|\\
&=\sup_{x\in\R^d,|u|=1}|\E[\nabla h(X_t^{x})\nabla_{u} X^{x}_t]|\\
&\leq \| \nabla h \|e^{-mt},
\end{align*}
where the last inequality followed from Lemma \ref{jm}.

By direct calculations, for any $u_1,u_2\in\R^d$, we have
\begin{align*}
&\nabla_{u_2}\nabla_{u_1}P_th(x)\\
=&\E[\nabla^2h(X_t^{x})\nabla_{u_2} X^{x}_t\nabla_{u_1} X^{x}_t+\nabla h(X_t^{x})\nabla_{u_2}\nabla_{u_1} X^{x}_t].
\end{align*}

Similarly, by  Lemma \ref{jm} and the  Cauchy-Schwarz inequality, we have
\begin{align*}
\|\nabla^2 P_th\|\leq C(\| \nabla h \|+\| \nabla^2 h \|) e^{-\frac{m}{2}t}.
\end{align*}
The proof is complete.
\end{proof}

\begin{proof}[Proof of Proposition \ref{prop:onestep}]
For any $f\in \mathcal{C}^2(\R^d,\R)$ and $y,z\in\R^d$, by the second order Taylor formula, we have
\begin{align*}
f(z)-f(y)=\langle\nabla f(y),z-y\rangle+\int_0^1(1-u)\nabla^2 f(uz+(1-u)y)(z-y)^{\otimes 2}\dif u.
\end{align*}
For a given $x\in\R^d$, it follows that
\begin{align*}
f(z)-f(y)=&\langle\nabla f(x),z-y\rangle+\langle\nabla f(y)-\nabla f(x),z-y\rangle\\
&+\int_0^1(1-u)\nabla^2 f(uz+(1-u)y)(z-y)^{\otimes 2}\dif u\\
=&\langle\nabla f(x),z-y\rangle+\int_0^1\langle \nabla^2f(x+u(y-x))(y-x),z-y\rangle\dif u\\
&+\int_0^1(1-u)\nabla^2 f(uz+(1-u)y)(z-y)^{\otimes 2}\dif u
\end{align*}
Applying this expansion with $z=\tl{Y}_{\gamma}^{x}$ and $y=X_{\gamma}^{x}$  yields
\begin{align}\label{I0}
|\E f(\tl{Y}_{\gamma}^{x})-\E f(X_{\gamma}^{x})|=\mathcal{I}_1+\mathcal{I}_2+\mathcal{I}_3,
\end{align}
where
\begin{align*}
&\mathcal{I}_1=|\E[\langle\nabla f(x),\tl{Y}_{\gamma}^{x}-X_{\gamma}^{x}\rangle]|,\\
&\mathcal{I}_2=\left|\E\left[\int_0^1\langle \nabla^2f(x+u(X_{\gamma}^{x}-x))(X_{\gamma}^{x}-x),\tl{Y}_{\gamma}^{x}-X_{\gamma}^{x}\rangle\dif u\right]\right|,\\
&\mathcal{I}_3=\left|\E\left[\int_0^1(1-u)\nabla^2 f(u\tl{Y}_{\gamma}^{x}+(1-u)X_{\gamma}^{x})(\tl{Y}_{\gamma}^{x}-X_{\gamma}^{x})^{\otimes 2}\dif u\right]\right|.
\end{align*}

It follows from Lemma \ref{one} that
\begin{align}
\mathcal{I}_2\leq&\|\nabla^2 f\|\|X_{\gamma}^{x}-x\|_2\|X_{\gamma}^{x}-\tilde{Y}_{\gamma}^{x}\|_2\leq C(1+|x|^2)\|\nabla^2 f\|\gamma^{2},\label{1345sf1}\\
\mathcal{I}_3\leq&C\|\nabla^2 f\|\|X_{\gamma}^{x}-\tilde{Y}_{\gamma}^{x}\|_2^2\leq C(1+|x|^2)\|\nabla^2 f\|\gamma^{3}.\label{1345sf2}
\end{align}

Next we estimate term $\mathcal{I}_1$.
Note that 
\begin{align}\label{1413sf}
|\E[\tl Y_{\gamma}^x-X_{\gamma}^x]|\leq\left|\int_0^\gamma\E[\nabla U(X_s)-\nabla U(x)]\dif s\right|
+\frac{\gamma}{K}\sum_{i=0}^{K-1}|\E[\nabla U(\hat{Y}_i)-\nabla U(x)]|,
\end{align}
where
\begin{align*}
\hat{Y}_i=x-\frac{i\gamma}{K}\nabla U(x)+\sqrt{\frac{2i\gamma}{K}}\xi_{i}, \quad i=0,1,\ldots,K-1.
\end{align*}

For the first term, using It\^{o}'s formula, we have for all $s\in[0,\gamma)$,
\begin{align}\label{ito}
\nabla U(X_s)-\nabla U(x)=&\int_0^s\{-\nabla^2U(X_r)\nabla U(X_r)+\vec{\Delta}(\nabla U)(X_r)\}\dif r\nonumber\\
&+\sqrt2\int_0^s\nabla^2U(X_r)\dif B_r.
\end{align}

By \eqref{1352sf},  \eqref{Xmoment} and for all $t\geq0$, $1-e^{-t}\leq t$, we have 
\begin{align}\label{1436sf}
\left|\int_0^\gamma\E[\nabla U(X_s)-\nabla U(x)]\dif s\right|=&\left|\int_0^\gamma\int_0^s\E\{-\nabla^2U(X_r)\nabla U(X_r)+\vec{\Delta}(\nabla U)(X_r)\}\dif r\dif s\right|\nonumber\\
\leq&\int_0^\gamma\int_0^s\{L^2\E|X_r|+d^2\tl L^2\}\dif r\dif s\nonumber\\
\leq&\int_0^\gamma\int_0^s\{L^2\sqrt{|x|^2+2dr}+d^2\tl L^2\}\dif r\dif s\nonumber\\
\leq&\int_0^\gamma\int_0^s\{L^2(|x|+\sqrt{2d}r^{1/2})+d^2\tl L^2\}\dif r\dif s\nonumber\\
=&\frac{1}{2}(L^2|x|+d^2\tl L^2)\gamma^2+\frac{4\sqrt2}{15}\sqrt dL^2\gamma^{\frac{5}{2}}.
\end{align}

To estimate the second term in RHS of \eqref{1413sf}, we construct an auxiliary process
\begin{align*}
\bar Y_{i,t}=x-t\nabla U(x)+\sqrt 2B_t, \quad t\in[0,\frac{i\gamma}{K}].
\end{align*}
Then $\bar Y_{i,\frac{i\gamma}{K}}\overset{d}{=}\hat Y_i$ and we have $\E[\nabla U(\hat Y_i)-\nabla U(x)]=\E[\nabla U(\bar Y_{i,\frac{i\gamma}{K}})-\nabla U(x)]$.

From It\^{o}'s formula,
\begin{align*}
|\E[\nabla U(\bar Y_{i,\frac{i\gamma}{K}})-\nabla U(x)]|
\leq&\left|\int_0^{\frac{i\gamma}{K}}\E[-\nabla^2U(\bar Y_{i,s})\nabla U(x)+\vec{\Delta}(\nabla U)(\bar Y_{i,s})]\dif s\right|\\
\leq&\int_0^{\frac{i\gamma}{K}} L^2|x|+d^2\tl L^2\dif s\\
\leq&(L^2|x|+d^2\tl L^2)\frac{i\gamma}{K},
\end{align*}
which implies 
\begin{align}\label{1437sf}
\frac{\gamma}{K}\sum_{i=0}^{K-1}|\E[\nabla U(\hat{Y}_i)-\nabla U(x)]\leq&\frac{\gamma}{K}\sum_{i=0}^{K-1}(L^2|x|+d^2\tl L^2)\frac{i\gamma}{K}\nonumber\\
=&\frac{K-1}{2K}(L^2|x|+d^2\tl L^2)\gamma^2.
\end{align}
Combining \eqref{1413sf}, \eqref{1436sf} and \eqref{1437sf}, we have   that
\begin{align*}
\mathcal{I}_1=|\E[\langle\nabla f(x),\tl{Y}_{\gamma}^{x}-X_{\gamma}^{x}\rangle]|\leq \|\nabla f\|
\left\{(1-\frac{1}{2K})(L^2|x|+d^2\tl L^2)\gamma^2+\frac{4\sqrt2}{15}\sqrt dL^2\gamma^{\frac{5}{2}}\right\}
\end{align*}

Collecting the estimates of $\mathcal{I}_1, \ \mathcal{I}_2$ and $\mathcal{I}_3$, we get the desired result. 

The proof is complete.

\end{proof}

\end{appendix}

\section*{Acknowledgements}
This work was partly  supported by National Natural Science Foundation of China (Grant Nos. 12271475 and U23A2064).

\end{document}